\documentclass[10pt]{article}
\usepackage{amsmath,amssymb,amsthm,amscd}
\numberwithin{equation}{section}

\def\Bl{\biggl(}
\def\Br{\biggr)}
\def\p{\partial}
\def\b{\bar}
\def\a{\alpha}
\def\g{\gamma}

\def\e{\eta}
\def\th{\theta}

\def\o{\omega}
\def\l{\lambda}

\def\tr{\rm tr}

\def\g{\mathfrak g}

\def\RR{{\mathbb R}}
\def\CC{{\mathbb C}}

\newtheorem{prop}{Proposition}[section]
\newtheorem{theo}[prop]{Theorem}
\newtheorem{lem}[prop]{Lemma}
\newtheorem{cor}[prop]{Corollary}
\newtheorem{rem}[prop]{Remark}

\newtheorem{defi}[prop]{Definition}

\newtheorem{q}[prop]{Question}
\def\begeq{\begin{equation}}
\def\endeq{\end{equation}}

\def\and{\quad{\rm and}\quad}

\let\lra=\longrightarrow

\def\mapright\#1{\,\smash{\mathop{\lra}\limits^{\#1}}\,}

\begin {document}

\bibliographystyle{plain}
\title{Ricci flow on K\"ahler-Einstein surfaces}
\author{X. X. Chen and G. Tian}
\date{July 26, 2000}
\maketitle

 \tableofcontents
\section{Introduction}

In the last two decades, the Ricci flow, introduced by R.
Hamilton in \cite{Hamilton82}, has been a subject of intense
study. The Ricci flow provides an indispensable tool of deforming
Riemannian metrics towards canonical metrics, such as Einstein
ones. It is hoped that by deforming a metric to a canonical
metric, one can further understand geometric and topological
structures of underlying manifolds. For instance, it was proved
\cite{Hamilton82} that any closed 3-manifold of positive Ricci
curvature is diffeomorphic to a spherical space form. We refer
the readers to \cite{Hamilton93} for more information.

If the underlying manifold is a K\"ahler manifold, the Ricci flow
preserves the K\"ahler class. It follows that the Ricci flow can
be reduced to a fully nonlinear parabolic equation on functions
(cf. Section 2 for details). Usually, this reduced flow is called
the K\"ahler Ricci flow. Unlike the Ricci flow in the real case,
it  can be  proved directly that
 the K\"ahler Ricci flow always has a global solution (cf. \cite{Cao85}).
 Following a similar calculation of Yau \cite{Yau78}, Cao \cite{Cao85} proved
that the solution converges to a K\"ahler-Einstein metric if the
first Chern class of the underlying K\"ahler manifold is zero or
negative. Consequently, he re-proved the famous Calabi-Yau
theorem\cite{Yau78}. On the other hand, if the first Chern class
of the underlying K\"ahler manifold is positive, the solution of a
K\"ahler Ricci flow may not converge to any K\"ahler-Einstein
metric. This is because there are compact K\"ahler manifolds with
positive first Chern class which do not admit any
K\"ahler-Einstein metrics (cf. \cite{futaki83} \cite{Tian97}). A
natural and challenging problem is whether or not the K\"ahler
Ricci flow on a compact K\"ahler-Einstein manifold converges to a
K\"ahler-Einstein metric. It was proved by S. Bando
\cite{Bando84} for 3-dimensional K\"ahler manifolds and by N. Mok
\cite{Mok88} for higher dimensional K\"ahler manifolds that the
positivity of bisectional curvature is preserved under the
K\"ahler Ricci flow. A long standing problem in the study of the
Ricci flow is whether or not the K\"ahler Ricci flow converges to
a K\"ahler-Einstein metric if the initial metric has positive
bisectional curvature? In view of the solution of the Frankel
conjecture by S. Mori \cite{Mori79} and Siu-Yau \cite{Siuy80}, we
suffice to study this problem on a K\"ahler manifold which is
biholomorphic to $\CC P^n$. Since $\CC P^n$ admits a
K\"ahler-Einstein metric, the above problem can be restated as
follows: on a compact K\"ahler-Einstein manifold, does the
K\"ahler Ricci flow converge to a K\"ahler-Einstein metric? This
problem was completely solved by R. Hamilton in the case of
Riemann surfaces (cf. \cite{Hamilton88}). We also refer the
readers to B. Chow's papers \cite{Chow91} for more developments
on this problem. In this paper, we give an affirmative answer to
this problem in dimension two.

\begin{theo} \footnote{In a subsequent paper \cite{chentian002},
 we will prove the same theorem for all dimensions. The proof for higher
dimensions needs new ingredients. Both results were announced in
\cite{chentian000}.} Let $M$ be a K\"ahler-Einstein surface with
positive scalar curvature. If the initial metric has nonnegative
bisectional curvature and positive at least at one point, then
the K\"ahler Ricci flow will converge exponentially fast to a
K\"ahler-Einstein metric with constant bisectional curvature.
\end{theo}

\begin{cor}
The space of K\"ahler metrics with non-negative bisectional
curvature (and positive at least at one point) is path-connected.
Moreover, the space of metrics with non-negative curvature
operator (and positive at least at one point)  is also
path-connected.
\end{cor}

\begin{rem}
Using the same arguments, we can also prove the version of
our main theorem for K\"ahler orbifolds.
\end{rem}

\begin{rem}
What we really need is that the Ricci curvature is
positive. Since the condition on Ricci may not be preserved
under the Ricci flow, in order to have the positivity of
the Ricci curvature, we will use the fact that
the positivity of the bisectional curvature is
preserved.
\end{rem}

\begin{rem} We need the assumption on the existence of K\"ahler-Einstein
metric because we will use a nonlinear inequality from
\cite{tian98}. Such an inequality is nothing but the
Moser-Trudinger-Onofri type if the K\"ahler-Einstein manifold is
the Riemann sphere.
\end{rem}

The typical method in studying the Ricci flow depends on
pointwise bounds of the curvature tensor by using its evolution
equation as well as the blow-up analysis. In order to prevent
formation of singularities, one blows up the solution of the Ricci
flow to obtain profiles of singular solutions. Those profiles
involve Ricci solitons and possibly more complicated singular
models. Then one tries to exclude formation of singularities by
checking that these solitons or models do not exist under
appropriate global geometric conditions. It is a common sense
that it is very difficult to detect how the global geometry
affects those singular models even for a very simple manifold
like $\CC P^2$. The first step is to classify those singular
models and hope to find their geometric information. Of course,
it is already a very big task. There have been many exciting
works on these (cf. \cite{Hamilton93}).

Our new contribution is to find a set of new functionals which are the Lagrangians
of certain new curvature equations involving various
symmetric functions of the Ricci curvature.
We show that these functionals decrease essentially along the
K\"ahler Ricci flow and have uniform lower bound.
By computing their derivatives, we can obtain certain integral bounds on
curvature of metrics along the flow.

For the readers' convenience, we will discuss more on these new
functionals. Let $M$ be a compact K\"ahler manifold with positive
first Chern class $c_1(M)$ and $\omega$ be a fixed K\"ahler
metric on $M$ with the K\"ahler class $c_1(M)$. Consider the
following expansion
 \begin{equation}
\left(\omega + t Ric(\omega)\right)^n = \left(\displaystyle
\sum_{k=0}^n\; \sigma_k(\omega) t^k \right)\omega^n,
 \label{symm}
 \end{equation}
where $\sigma_k(\omega)$ is the $k-$th symmetric polynomial of the
Ricci tensor ${\rm Ric}(\omega).\;$
  Then we say that a K\"ahler metric $\omega$ is of
extremal k-th symmetric Ricci curvature ($k=0,1, \cdots, n$) if
$\sigma_k(\omega)$ satisfies
\begin{equation}
\Delta \sigma_{k} (\omega) - \frac{n-k}{k+1} \sigma_{k+1}(\omega)
= c_k,
\end{equation}
where $c_k$ is a constant determined by $c_1(M)$ and the K\"ahler
class $[\omega]$. Clearly, the extremal 0-symmetry means constant
scalar curvature. When the first Chern class $c_1(M)$ of $M$ is
positive and $\omega$ represents $c_1(M)$, a K\"ahler metric with
constant scalar curvature is of constant Ricci curvature and
consequently, has extremal k-th symmetric Ricci curvature for all
$k$. In general, a K\"ahler metric of constant scalar curvature
may not have extremal k-th symmetric Ricci curvature for $k> 1$.

Our new functionals $E_k$ are simply the Lagrangians of the above
Ricci curvature equations (cf. section 4 for details). When $k=0$,
the functional $E_0$ is nothing but the K-energy of T. Mabuchi. We
will prove that the derivative of each $E_k$ along an orbit of
automorphisms gives rise to a holomorphic invariant $\Im_k$,
including the well-known Futaki invariant as a special one. When
$M$ admits a K\"ahler-Einstein metric, all these invariants
$\Im_k$ vanish, so the functionals $E_k$ are invariant under the
action of automorphisms.

Next we will prove that these $E_k$ are bounded from below. This
can be achieved by making use of a fully nonlinear inequality
from \cite{tian98} (cf. Section 5).  But in order to apply this
inequality, we have to adjust the fixed K\"ahler-Einstein metric
 so that the evolved K\"ahler metrics are centrally positioned with
respect to the adjusted K\"ahler-Einstein metrics, that is, the
K\"ahler  potentials between the two evolved metrics are
orthogonal to the first eigenspace of the evolved
K\"ahler-Einstein metrics (cf. Section 6). It causes some extra
difficulties in the proof of our main theorem (particularly in
higher dimensions).

 Next we will compute the derivatives of
$E_k$ along the K\"ahler-Ricci flow. Recall that the K\"ahler
Ricci flow is given by
\begin{equation}
   {{\partial \varphi} \over {\partial t }} =  \log {(\omega + \sqrt{-1}\partial \overline\partial
\varphi)^n \over {\omega}^n } + \varphi - h_{\omega},
\end{equation}
where $h_\omega$ depends only $\omega$. The derivatives of these
functionals are all bounded uniformly from above along the
K\"ahler Ricci flow. Furthermore, we found that $E_0$ and $E_1$
decrease along the K\"ahler Ricci flow. These play a very
important role in this and the subsequent paper. We can derive
from these properties of $E_k$ integral bounds on curvature, e.g.
for almost all K\"ahler metrics $\omega_{\varphi(t)}$ along the
flow, we have
\begin{equation}
  \int_M\; R({\omega_{\varphi(t)}})\; {\rm Ric({\omega_{\varphi(t)}})}^k
   \wedge{\omega_{\varphi(t)}}^{n-k} \le C, \;k=1, \cdots, n,
\label{eq:keyestimate1}
\end{equation}
and
\begin{equation}
 \int_M\; (R(\omega_{\varphi(t)})-r)^2 \;{\omega_{\varphi(t)}}^n \rightarrow 0,
\label{eq:keyestimate2}
\end{equation}
where $R({\omega_{\varphi(t)}})$ denotes the scalar curvature and
$r$ is the average scalar curvature.

In principle, one can then follow Hamilton's arguments in the
case of Riemann surfaces. But we need to do some changes since
the sectional curvature may not be positive and we can not apply
Klingenberg's estimate on injectivity radius. We will generalize
Klingenberg's estimate to K\"ahler manifolds of positive
bisectional curvature. Then, combining the above integral bounds
on the curvature with Cao's Harnack inequality and the
generalization of Klingenberg's estimate, we can bound the
curvature uniformly along the K\"ahler Ricci flow in the case of
K\"ahler-Einstein surfaces. Then it is quite routine to prove the
convergence to the K\"ahler-Einstein metric. For higher
dimensions, one has to develop new techniques in order to get the
curvature bound. We will do it in a subsequent paper
\cite{chentian002}.

The organization of our paper is roughly as follows: In Section 2,
we review briefly some basics in K\"ahler geometry and necessary
information on the K\"ahler Ricci flow. In Section 3, we discuss
two important energy functionals. In Section 4, we introduce a set
of new functionals as we have briefly described in the above. In
Section 5, we prove that these functionals are invariant on any
K\"ahler-Einstein manifolds.  In Section 6, we modify the evolved
K\"ahler metrics to obtain desired integral estimates on the
curvature . In Section 7,8,9, we will bound the scalar curvature
uniformly along the K\"ahler Ricci flow. In Section 10, we prove
the exponentially convergence. In Section 11, we make some
concluding remarks and propose some open questions.

\section{Basic K\"ahler Geometry}
\subsection{Notations in K\"ahler geometry}
Let $M$ be an $n$-dimensional compact K\"ahler manifold. A K\"ahler metric can be given by its
K\"ahler form $\omega$ on $M$. In local coordinates $z_1, \cdots, z_n$, this $\omega$ is of
the form
\[
\omega = \sqrt{-1} \displaystyle \sum_{i,j=1}^n\;g_{i \overline{j}} d\,z^i\wedge d\,z^{\overline{j}},
\]
where $\{g_{i\overline {j}}\}$ is a positive definite Hermitian matrix function.
The K\"ahler condition requires that $\omega$ is a closed positive
(1,1)-form. In other words, the following holds
\[
 {{\partial g_{i \overline{k}}} \over
{\partial z^{j}}} =  {{\partial g_{j \overline{k}}} \over
{\partial z^{i}}}\qquad {\rm and}\qquad {{\partial g_{k
\overline{i}}} \over {\partial z^{\overline{j}}}} = {{\partial
g_{k \overline{j}}} \over {\partial
z^{\overline{i}}}}\qquad\forall\;i,j,k=1,2,\cdots, n.
\]
The K\"ahler metric corresponding to $\omega$ is given by
\[
 \sqrt{-1} \;\displaystyle \sum_1^n \; {g}_{\alpha \overline{\beta}} \;
d\,z^{\alpha}\;\otimes d\, z^{ \overline{\beta}}.
\]
For simplicity, in the following, we will often denote by $\omega$ the corresponding K\"ahler metric.
The K\"ahler class of $\omega$ is its cohomology class $[\omega]$ in $H^2(M,\RR).\;$
By the Hodge theorem, any other K\"ahler
metric in the same K\"ahler class is of the form
\[
\omega_{\varphi} = \omega + \sqrt{-1} \displaystyle
\sum_{i,j=1}^n\; {{\partial^2 \varphi}\over {\partial z^i
\partial z^{\overline{j}}}} \;d\,z_i \wedge d\,z_{\bar j}
> 0
\]
for some real value function $\varphi$ on $M.\;$ The functional
space in which we are interested (often referred as the space
of K\"ahler potentials) is
\[
{\cal P}(M,\omega) = \{ \varphi \;\mid\; \omega_{\varphi} = \omega
+ \sqrt{-1}
 {\partial} \overline{\partial} \varphi > 0\;\;{\rm on}\; M\}.
\]
Given a K\"ahler metric $\omega$, its volume form  is
\[
  \omega^n = \left(\sqrt{-1} \right)^n \det\left(g_{i \overline{j}}\right)
 d\,z^1 \wedge d\,z^{\overline{1}}\wedge \cdots \wedge d\,z^n \wedge d\,z^{\overline{n}}.
\]
Its Christoffel symbols are given by
\[
  \Gamma^k_{i\,j} = \displaystyle \sum_{l=1}^n\;g^{k\overline{l}} {{\partial g_{i \overline{l}}} \over
{\partial z^{j}}} ~~~{\rm and}~~~
\Gamma^{\overline{k}}_{\overline{i}\,\overline{j}} =
\displaystyle \sum_{l=1}^n\;g^{\overline{k}l} {{\partial g_{l
\overline{i}}} \over {\partial z^{\overline{j}}}},
\qquad\forall\;i,j,k=1,2,\cdots, n.
\]
The bisectional curvature tensor is
\[
 R_{i \overline{j} k \overline{l}} = - {{\partial^2 g_{i \overline{j}}} \over
{\partial z^{k} \partial z^{\overline{l}}}} + \displaystyle
\sum_{p,q=1}^n g^{p\overline{q}} {{\partial g_{i \overline{q}}}
\over {\partial z^{k}}}  {{\partial g_{p \overline{j}}} \over
{\partial z^{\overline{l}}}}, \qquad\forall\;i,j,k,l=1,2,\cdots,
n.
\]
We say that $\omega$ is of nonnegative bisectional curvature if
\[
 R_{i \overline{j} k \overline{l}} v^j v^{\overline{i}} w^l w^{\overline{k}}\geq 0
\]
for all non-zero vectors $v$ and $w$ in the holomorphic tangent
bundle of $M$. The bisectional curvature and the curvature tensor
can be mutually determined by each other (cf. The appendix for
more information). The Ricci curvature of $\omega$ is locally
given by
\[
  R_{i \overline{j}} = - {{\partial}^2 \log \det (g_{k \overline{l}}) \over
{\partial z_i \partial \bar z_j }} .
\]
So its Ricci curvature form is
\[
  {\rm Ric}(\omega) = \sqrt{-1} \displaystyle \sum_{i,j=1}^n \;R_{i \overline{j}}(\omega)
d\,z^i\wedge d\,z^{\overline{j}} = -\sqrt{-1} \partial \overline{\partial} \log \;\det (g_{k \overline{l}}).
\]
It is a real, closed (1,1)-form. Recall that $[\omega]$ is a
canonical K\"ahler class if this Ricci form is cohomologous to
$\lambda \;\omega,\; $ for some constant $\lambda$

\subsection{ The K\"ahler Ricci flow }

    Now we assume that the first Chern class $c_1(M)$ is positive.
The Ricci flow (see for instance \cite{Hamilton82} and
\cite{Hamilton86}) on a K\"ahler manifold $M$ is of the form
\begin{equation}
  {{\partial g_{i \overline{j}}} \over {\partial t }} = g_{i \overline{j}}
  - R_{i \overline{j}}, \qquad\forall\; i,\; j= 1,2,\cdots ,n.
\label{eq:kahlerricciflow}
\end{equation}
If we choose the initial K\"ahler metric $\omega$ with $c_1(M)$
as its K\"ahler class. Then the flow (2.1) preserves the K\"ahler
class $[\omega]$. It follows that on the level of K\"ahler
potentials, the Ricci flow becomes
\begin{equation}
   {{\partial \varphi} \over {\partial t }} =  \log {{\omega_{\varphi}}^n \over {\omega}^n } + \varphi - h_{\omega} ,
\label{eq:flowpotential}
\end{equation}
where $h_{\omega}$ is defined by
\[
  {\rm Ric}(\omega)- \omega = \sqrt{-1} \partial \overline{\partial} h_{\omega}, \; {\rm and}\;\displaystyle \int_M\;
  (e^{h_{\omega}} - 1)  {\omega}^n = 0.
\]
As usual, the flow (2.2) is referred as the K\"ahler Ricci flow
on $M$. Differentiating on both sides of equation
(\ref{eq:flowpotential}) on $t$, we obtain

\[
 {{\partial } \over {\partial t }}\; {{\partial \varphi} \over {\partial t }} =
 \triangle_{\varphi}  {{\partial \varphi} \over {\partial t }} +  {{\partial \varphi} \over {\partial t }},
\]
where $\triangle_{\varphi}$ is the Laplacian operator of the metric $\omega_{\varphi}.\;$
Then it follows from the standard Maximum Principle
\begin{lem} Along the K\"ahler Ricci flow (\ref{eq:kahlerricciflow}),
$| {{\partial \varphi} \over {\partial t }}|$ grows at most exponentially.
\end{lem}
In particular, the $C^0$-norm of $\varphi$ can be bounded by a
constant depending only $t$. Using this fact and following Yau's
calculation in \cite{Yau78}, one can prove   that for any initial
metric with K\"ahler class $c_1(M)$, the evolution equation (2.2)
has a global solution for all time $0\leq t< \infty$ (cf.
\cite{Cao85}).

\subsection{Preservation of nonnegative bisectional curvature}

The K\"ahler Ricci flow induces an evolution equation on the
bisectional curvature
\[
\begin{array}{lcl} {{\partial }\over {\partial t}} R_{i \overline{j} k \overline{l}}
& = & \bigtriangleup R_{i \overline{j} k \overline{l}} + R_{i \overline{j} p \overline{q}} R_{q \overline{p} k
\overline{l}} - R_{i \overline{p} k \overline{q}}
R_{p \overline{j} q \overline{l}} + R_{i \overline{l} p \overline{q}} R_{q \overline{p} k \overline{j}} +
R_{i \overline{j} k \overline{l}}\\
& &\qquad  -{1\over 2} \left( R_{i \overline{p}}R_{p \overline{j} k \overline{l}} +
R_{p \overline{j}}R_{i \overline{p} k \overline{l}} + R_{k \overline{p}}R_{i \overline{j} p \overline{l}} +
R_{p \overline{l}}R_{i \overline{j} k \overline{p}} \right).
\end{array}
\]
Similarly, we have evolution equation for the Ricci tensor and
the scalar curvature
\[
{{\partial }\over {\partial t}} R_{i \overline{j}} =  \bigtriangleup R_{i \overline{j}}
+ R_{l \overline{k}}R_{i \overline{j} k \overline{l}} - R_{i \overline{k}}R_{ k \overline{j}},
\]
and
\[{{\partial }\over {\partial t}} R =  \bigtriangleup R
+ |{\rm Ric}|^2- R.\]

The following theorem was proved by S. Bando for 3-dimensional
compact K\"ahler manifolds. This  was late by N. Mok in
\cite{Mok88} for all K\"ahler manifolds. Their proof used
Hamilton's Maximum Principle for tensors. The proof for higher
dimensions is quite intrigue.

\begin{theo} \cite{Bando84}  \cite{Mok88} Under the K\"ahler Ricci flow, if the initial metric
has nonnegative bisectional curvature, then the evolved metrics
also have non-negative bisectional curvature. Furthermore, if the
bisectional curvature of the initial metric is positive at least
at one point, then the evolved metric has positive bisectional
curvature at all points.
\end{theo}

Previously, R. Hamilton proved (by using his Maximum principle
for tensors)

\begin{theo} Under the K\"ahler Ricci flow, if the initial metric
has nonnegative curvature operator, then the evolved metrics also
have non-negative curvature operator. Furthermore, if the
curvature operator of the initial metric is positive at least at
one point, then the evolved metric has positive curvature
operator at all points.
\end{theo}

It is still interesting to see if similar conclusion holds for
sectional curvature, that is, if the initial metric has
nonnegative sectional curvature, do evolved metrics along the
K\"ahler Ricci flow have nonnegative sectional curvature? If so,
our theorem will imply that there is no exotic K\"ahler metric
with positive sectional curvature
on complex projective spaces.\\

\section{Generalized energy functionals}
In this section, we will introduce some generalized energy
functionals $J_{\omega}$, $F_{\omega}$ and $\nu_{\omega}$. The
second functional was first used in \cite{DingTian93}, while the
3rd one was introduced by T. Mabuchi. These are all useful
functionals in K\"ahler geometry. We then review some known
properties of $F_{\omega}$ and $\nu_{\omega}$, such as a) they
both decrease under the K\"ahler Ricci flow; b) They are both
invariant under automorphisms on any K\"ahler-Einstein manifolds.

\subsection{A nonlinear inequality}

Recall that the generalized energy:
\begin{equation}
J_{\omega}(\varphi) = {1 \over V} \displaystyle
\sum_{i=0}^{n-1}\;\displaystyle \int_M\;{{i+1}\over{n+1}}\;
\sqrt{-1}\;\partial \varphi \wedge \overline{\partial}\varphi
\wedge \omega^i \wedge {\omega_{\varphi}}^{n-1-i}.
\label{eq:Jdefinition}
\end{equation}
where $V = \displaystyle \int_M\; \omega^n = [\omega]^n ([V]) $ and
$\omega_{\varphi} = \omega + \sqrt{-1} \partial  \overline{\partial}
\varphi.\;$ This is clearly a positive functional.
When $n=1$, it is just the standard Dirichlet energy
\[
J_\o (\varphi ) = {1 \over 2V} \int _M \;\sqrt{-1}\;\p \varphi
\wedge \b \p \varphi = {1 \over 2V} \int _M |\p \varphi | ^2 \o.
\]

If $n=2,\;$ we have
\[
J_\o (\varphi ) = {1 \over 3V} \int _M \;\sqrt{-1}\;\p \varphi
\wedge \b \p \varphi \wedge \o _\varphi + {2 \over 3V} \int _M
\;\sqrt{-1}\;\p \varphi \wedge \b \p \varphi \wedge \o.
\]

Taking derivative of $J_\omega$ along a path $\varphi(t) \in
{\cal P}(M,\omega) $, we arrive at
\[
{{d\, J_{\omega}(\varphi)}\over { d\,t}} = - {1 \over V}
\;\displaystyle \int_M \; {{\partial \varphi}\over { d\,t}}
\left({\omega_{\varphi}}^n - {\omega}^n \right).
\]

Alternatively, one can use this formula to define $J_\omega$. From
this formula, one can see that $J_{\omega}$ does not satisfy the
cocycle condition. Recall that the functional $F_\o$ is defined by
 \[
F_{\omega} (\varphi) = J_{\omega}(\varphi) - {1 \over V}
\;\displaystyle \int_M \; \varphi\;\; {\omega}^n - \log \left(
{1\over V} \;\displaystyle \int_M \; e^{h_{\omega} -
\varphi}\right).
\]

It satisfies the cocycle condition and its critical points are
K\"ahler-Einstein metrics. If $n=1$, then $M = S^2$ and
\[
F_\o (\varphi ) ={1 \over 2V} \int _{S^2} |\p \varphi |^2 - {1
\over V} \int _{S^2}
 \varphi \o - \log {1 \over V} \int _{S^2} e^{h_\o -\varphi } \o .
\]

This is precisely the functional in studying L. Nirenbberg's
problem of prescribing the Gauss curvature on $S^2$.

Suppose that $M$ has positive first Chern class and admits a
K\"ahler-Einstein metric. Then there is a K\"ahler-Einstein metric
$\omega_1$ such that ${\rm Ric}(\omega_1) = \omega_1$. We will
denote by $\Lambda_1$ the space of eigenfunctions with eigenvalue
one if one is an eigenvalue of the K\"ahler-Einstein metric
$\omega_1$. If one is not an eigenvalue, we simply put
$\Lambda_1$ to be $\{0\}$. By $\phi \perp \Lambda_1$, we mean
$\int_M\; \phi \psi {\omega_1}^n = 0$ for all $\psi \in
\Lambda_1$. Note that if $M$ admits no holomorphic vector fields,
then $\Lambda_1=\{0\}$ and $\phi \perp \Lambda_1 $ is
automatically true.

The following inequality plays an important role in our proof.

\begin{theo}(Tian) \cite{Tian97} Let $M$ be given as above and $\omega_1$ be any K\"ahler metric
with $c_1(M)$ as its K\"ahler class. Then there exist constants
$\delta = \delta(n) $ and $c = c(M, \omega_1)\geq 0$ such that for
any $\phi \in \cal H$ which satisfies $\phi \perp \Lambda_1,$ we
have
\[
  F_{\omega_1} (\phi) \geq J_{\omega_1}(\phi)^{\delta} - c
\]
which is the same as
\[
  {1\over V} \displaystyle \int_M e^{-\phi} {\omega_1}^n \leq C
e^{J_{\omega_1}(\phi) - {1\over V} \int_M \phi {\omega_1}^n -
{J_{\omega_1}(\phi)}^{\delta}}.
\]
\end{theo}

\begin{rem} This inequality was first proved under an extra condition, which
was removed later in \cite{TianZhu001}.
\end{rem}

\begin{rem} Since the difference of $J_\omega$ and $J_{\omega_1}$ (resp. $F_\omega$ and $F_{\omega_1}$)
is bounded by a constant depending only on $\omega$ and $\omega_1$, the inequality in the above
theorem holds irrelevant of choices of initial metrics.
\end{rem}

Inspired by the work of Donaldson \cite{Dona87}, T. Mabuchi
introduced the K-energy.

\begin{defi}(Mabuchi \cite{Ma87}) For any $\varphi (t) \in \cal P,$ the derivative of the K-energy
along this path $\varphi(t)$ is:
\[
{{d\,}\over {d\,t}} \nu_{\omega}(\varphi(t)) = - { 1 \over V}\;
\int_M\; {{\partial \varphi}\over {\partial t}} \;(R(\varphi(t))
- r) \;{\omega_{\varphi}}^n,
\]
where $r$ is the average value of the scalar curvature $r =
{{[c_1(M)]\cdot [\omega]^{n-1}} \over {[\omega]^{n}}}.\;$
\end{defi}

It was found in \cite{Tian94} that the K-energy can be expressed as
\begin{eqnarray}
\nu_{\omega}(\varphi) & = & \displaystyle {1 \over V}\; \int_M
\log\left ({\omega_\varphi^n \over \omega^n}\right
)\omega^n_\varphi -\displaystyle {1 \over V}\;\int_M h_\omega
(\omega^n -\omega_\varphi^n) \nonumber \\
& & \qquad  - \displaystyle {\sqrt{-1} \over V}\; \displaystyle
\sum_{i=0}^{n-1}\; \displaystyle \int_M \;\p \varphi \wedge \bar
\p \varphi \wedge \omega^i \wedge {\omega_{\varphi}}^{n-i-1}.
 \label{eq:kenergyexp0}
\end{eqnarray}

It was also observed in \cite{DingTian93} that
$$\nu_\omega (\varphi) \ge F_\omega(\varphi) -{1\over V}\;\displaystyle \int_M h_\omega \omega^n.$$
Combining this with the above theorem, we get

\begin{cor} \cite{Tian97} Suppose that $c_1(M) > 0$ and there is a K\"ahler-Einstein metric on $M$.
Then for any K\"ahler metric $\omega$ with $c_1(M)$ as its
K\"ahler class, there are constants $\delta = \delta(n) $ and $c
= c(M, \omega)\geq 0$ such that for any $\varphi \in {\cal
P}(M,\omega)$ which satisfies $\varphi \perp \Lambda_1,$ we have
\[
  \nu_{\omega} (\phi) \geq J_{\omega}(\phi)^{\delta} - c,
\]
\end{cor}

The following corollary will be crucial in our arguments.

\begin{cor}  Suppose that $c_1(M) > 0$ and there is a K\"ahler-Einstein metric on $M$.
Then for any function $\varphi \in {\cal P}(M,\omega)$
perpendicular to $\Lambda_1$, we have
$$\int_M \log\left ({\omega_\varphi^n \over \omega^n}\right )\omega^n_\varphi
\le C ( 1+ \nu_\omega (\varphi))^{1\over \delta},
$$
where $C$ is a constant depending only on $M$ and $\omega$.
\end{cor}

\subsection{Monotonicity  along the K\"ahler Ricci flow}

First we collect two simple facts which were known to experts in
the field for a while.

\begin{lem} Under the K\"ahler Ricci flow, $F_{\omega}$ decreases monotonely.
\end{lem}
\begin{proof} Let $c =  \log \left( {1\over V} \;\displaystyle \int_M \; e^{h_{\omega} - \varphi}\right).\;$ Then
\[
\begin{array}{lcl} {{d\,}\over {d\,t}} F_{\omega}(\varphi(t)) & = & -{1\over V}  \int_M\;{{\partial \varphi}\over {\partial t}} \;
 \left( \;{\omega_{\varphi}}^n  - e^{h_{\omega} - \varphi - c} \omega^n \right)\\
& = & - {1\over V} \int_M\;(\log {{\omega_{\varphi}}^n  \over
{{\omega}^n} } + \varphi - h_{\omega})\;
\left({{\omega_{\varphi}}^n  \over {\omega}^n}  - e^{h_{\omega} -
\varphi -c } \right) \;\omega^n \\  & = & - {1\over V}
\int_M\;(\log {{\omega_{\varphi}}^n  \over {{\omega}^n} } -
(h_{\omega}- \varphi -c ))\; \left({{\omega_{\varphi}}^n \over
{\omega}^n}  - e^{h_{\omega} - \varphi -c } \right) \;\omega^n\\
& \leq & 0.
\end{array}
\]
\end{proof}

Similarly, we have

\begin{lem} Under the K\"ahler Ricci flow, the K-energy $\nu_\o$ monotonely decreases!
\end{lem}
\begin{proof} By the definition, we have
\[
\begin{array}{lcl} {{d\,}\over {d\,t}} \nu_{\omega}(\varphi(t)) & = &
-{1\over V}
\int_M\;{{\partial \varphi}\over {\partial t}} \; (R(\omega_{\varphi(t)}) - r) \;{\omega_{\varphi_t}}^n\\
& = & {1\over V} \int_M\;(\log {{\omega_{\varphi_t}}^n  \over
{{\omega}^n} } + \varphi - h_{\omega})\;
(\bigtriangleup_{\varphi_t}(\log {{\omega_{\varphi_t}}^n  \over {\omega}^n}  + \varphi - h_{\omega}) ) \;{\omega_{\varphi}}^n \\
 & \leq & 0.
\end{array}
\]

The lemma follows.
\end{proof}

Next we want to prove that $F_\omega$ and
$\nu_{\omega}$ are both invariant under automorphisms on a
K\"ahler-Einstein manifold.

Recall that the Futaki invariant $f_M$ can be defined by
(see \cite{futaki83} )
$$
f_M (\o ,X) = \int _M X(h_\omega) \o  ^n,
$$
where $\omega$ is a K\"ahler metric with $c_1(M)$ as its K\"ahler
class and $X$ is a holomorphic vector field on $M$. Futaki proved
that the integral is independent of the choice of $\omega$, so it
gives rise to a holomorphic invariant. If $M$ admits a
K\"ahler-Einstein metric, then $f_M \equiv 0$.

Let $\Phi _t$ be a one-parameter group of automorphisms generated
by ${\rm Re} (X)$. Write $\omega_t=\Phi _t ^* \o = \o + \sqrt{-1}
\p \b \p \varphi _t$. We can further normalize $\varphi _t$ such
that $\int _M (e^{h_\o -\varphi _t} -1)\o ^n =0$. Then $h_{\o _t}
= \Phi _t ^* h_\o $. This implies that $\dot h_{\o _t} = {\rm
Re}(X) (h_{\o _t})$.

On the other hand, using the identity
$$
{\rm Ric} (\o _t) -\o _t = {\rm Ric} (\o) - \o -\sqrt{-1} \p \b \p
\log \Bl {{\o _t ^n} \over {\o ^n}} \Br -\sqrt{-1} \p \b \p
\varphi _t,
$$
we get
$$
h_{\o _t} =h_\o -\log \Bl {{\o _t ^n} \over {\o ^n}} \Br - \varphi
_t .
$$

Differentiating it with respect to $t$, we have
$$
\dot h_{\o _t} = -\Delta _{\o _t} \dot {\varphi _t} -  \dot
{\varphi _t}.
$$

Combining all these, we arrive at
$$
{d \over dt} F_\o (\varphi _t ) = {1\over V}\;{\rm Re}(f_M (X)).
$$

The following corollary is an immediate consequence of this.

\begin{lem}
The functional $F_\o$ is invariant under automorphisms if $f_M
\equiv 0$. In particular, it is true if $M$ is a K\"ahler-Einstein
manifold.
\end{lem}

Similarly, we have
\begin{lem} On a K\"ahler-Einstein manifold, $\nu_{\omega}$ is invariant under automorphisms.
\end{lem}

We deduce the following from the above

\begin{prop}  Suppose that $M$ admits a K\"ahler-Einstein metric.
Let $\varphi(t)$ ($t>0$) be a global solution of the K\"ahler Ricci flow and $\Psi(t)$
be a family of automorphisms of $M$. Write
$$\Psi_t^*\omega_{\varphi(t)} = \omega + \sqrt{-1}\p\b \p \psi(t).$$
Then $F_\o(\psi_t)$ and $\nu_{\omega}(\psi(t))$ are decreasing functions of $t$.
\end{prop}

Combining this with Tian's inequality last subsection, we get

\begin{cor}
Suppose that $\omega_1$ is a K\"ahler-Einstein metric with ${\rm Ric}(\omega_1)=\omega_1$.
Let $\varphi_t$ ($t>0$) be a global solution of the K\"ahler Ricci flow and $\Psi_t$
be a family of automorphisms of $M$. Write
$$\Psi_t^*\omega_{\varphi_t} = \omega + \sqrt{-1}\p\b \p \psi_t.$$
Let $\varphi(t)$ ($t>0$) be a global solution of the K\"ahler Ricci flow and $\Psi(t)$
be a family of automorphisms of $M$. Write
$$\Psi_t^*\omega_{\varphi(t)} = \omega + \sqrt{-1}\p\b \p \psi(t).$$
If $\psi(t)$ is perpendicular to the eigenspace of $\omega_1$ with eigenvalue one,
then
$$J_\o(\psi(t)) \le \nu_\o(\varphi(0)) + c,$$
where $c$ is a uniform constant.
\end{cor}

\subsection{Estimate on volume forms}

The following is the main result of this subsection.

\begin{prop} If ${\rm Ric}(\omega_{\varphi}) \geq 0,$ then there exists a uniform constant $C$ such that

\[
\displaystyle \inf_{M}\; \left( \log {\omega_{\varphi}^n \over
{\omega}^n } \right)(x)\geq - 4 C\, e^{2(1+ \displaystyle \int_M
\; \left( \log {\omega_{\varphi}^n \over {\omega}^n }
\right)\;{\omega_{\varphi}}^n )}.
 \]
\end{prop}
\begin{proof}
Choose any constant $c$ such that
\[
\displaystyle\;{1\over V} \int_{M}
\log {{\omega_{\varphi}}^n \over {\omega}^n } {\omega_{\varphi}}^n \leq c,\]
where $V=\int_M \omega^n$.

Put $\epsilon $ to be $e^{-2(1+c)}$.
Observe that
$$\left(\log {\omega_{\varphi}^n \over
{\omega}^n } \right) \omega_{\varphi}^n \geq - e^{-1}\omega^n,$$

we have
\[
\begin{array}{lcl}
c V & \geq & \displaystyle \int_{ \epsilon \omega_{\varphi_i}^n > \omega^n }  \left(\log
{\omega_{\varphi}^n \over {\omega}^n }\right)
\omega_{\varphi}^n
+ \displaystyle \int_{\epsilon \omega_{\varphi_i}^n \le \omega^n }
 \left(\log {\omega_{\varphi}^n \over {\omega}^n }\right) \omega_{\varphi}^n \\
& \geq & \displaystyle \int_{\epsilon \omega_{\varphi_i}^n > \omega^n}
\left( \log {1\over \epsilon} \right) \omega_{\varphi}^n +
 \displaystyle \int_{ \epsilon \omega_{\varphi_i}^n \le \omega^n }  (- e^{-1} \omega^n)\\
& > & 2 (1+c)\int_{\epsilon \omega_{\varphi_i}^n > \omega^n } \omega_\varphi^n - V.
\end{array}
\]

It follows that
$$\int_{\epsilon \omega_{\varphi_i}^n > \omega^n } \omega_\varphi^n < {V\over 2},$$
and consequently,
$$\int_{\omega^n \le 4 \omega_{\varphi_i}^n } \omega^n
\ge \epsilon \int_{{\epsilon \over 4}\omega^n \le
\epsilon \omega_{\varphi}^n \le \omega^n } \omega_\varphi^n > {\epsilon  V\over 4}.
$$

The Ricci curvature being non-negative implies that
\[
  \omega + \sqrt{-1}\; \partial \overline{\partial}
\left( h_{\omega}- \log {{\omega_{\varphi}}^n \over {\omega}^n }\right) \geq 0.
\]

Taking trace with respect to $\omega$, we get
\[
 n + \Delta \left( h_{\omega} - \log {{\omega_{\varphi}}^n \over {\omega}^n } \right) \geq  0,
\]
where $\Delta$ denotes the Laplacian of $\omega$.
Then by the Green formula, we have
\[
\begin{array}{cll}
&&\left( h_{\omega} - \log {\omega_{\varphi}^n \over {\omega}^n }\right)(x) \\
& = & {1\over V} \displaystyle
\int_M \left( h_{\omega} - \log {\omega_{\varphi}^n \over {\omega}^n }\right) \omega^n
- {1\over V} \displaystyle \int_M \Delta \left( h_{\omega} - \log {{\omega_{\varphi}}^n \over {\omega}^n }\right)
G(x,y) \omega^n (y)\\
& \leq &  {1\over V} \displaystyle
\int_M  \left( h_{\omega} - \log {{\omega_{\varphi}}^n \over {\omega}^n }\right) {\omega}^n   -  {n\over V} \displaystyle
\int_M G(x,y) \\
& \leq & {1\over V} \displaystyle
\int_M  \left( h_{\omega} - \log {{\omega_{\varphi}}^n \over {\omega}^n }\right) {\omega}^n   +  c',
\end{array}
\]
where $G(x,y) \geq 0 $ is a Green function of $\omega$. Note that we will always
denote by $c'$ a constant depending only on $\omega$ in this proof.

It follows from the above inequalities that
\[
\begin{array}{cll} & &
\displaystyle \inf_{M}\; \left( \log {\omega_{\varphi}^n \over
{\omega}^n } \right)\\ &\geq  &{1\over V} \displaystyle
\int_M  \left( \log {\omega_{\varphi}^n \over {\omega}^n } \right) {\omega}^n  - c'\\
&\geq& \inf_{M}\; \left( \log {\omega_{\varphi}^n \over {\omega}^n } \right)
{1\over V} \displaystyle \int_{\omega^n\ge 4 \omega_\varphi^n}
\omega^n - {\log 4 \over V} \displaystyle\int _{\omega^n \le 4\omega^n_\varphi} \omega^n - c'\\
& \geq & (1- {\epsilon \over 4}) \inf_{M}\; \left( \log
{\omega_{\varphi}^n \over {\omega}^n } \right) - c'. \end{array}
\]

Therefore, we have
\[
\displaystyle \inf_{M}\; \left( \log {\omega_{\varphi}^n \over {\omega}^n } \right)(x)\geq - 4 c' e^{2(1+c)}.
\]
By the way we choose the constant $c$ in the beginning of the
proof, we have
\[
\displaystyle \inf_{M}\; \left( \log {\omega_{\varphi}^n \over
{\omega}^n } \right)(x)\geq - 4 c' e^{2(1+ \displaystyle \int_M
\; \left( \log {\omega_{\varphi}^n \over {\omega_{\varphi}}^n }
\right)\;{\omega}^n )}.
 \]

The proposition is proved.
\end{proof}

\section{New functionals}
In this section, we introduce a family of new
functionals on the space of K\"ahler potentials ${\cal
P}(M,\omega)$. We will show that their derivatives
along the K\"ahler Ricci flow are bounded uniformly
from above.

\subsection{Definition of functionals $E_k$}

In this subsection, we introduce $E_k$ for $k = 0,1, \cdots, n$.

\begin{defi} For any $k=0,1,\cdots, n$, we define a functional $E_k^0$
on ${\cal P}(M,\omega)$ by
\[
  E_{k,\omega}^0 (\varphi) = {1\over V}\; \displaystyle \int_M\;  \left( \log {{\omega_{\varphi}}^n \over \omega^n}
   - h_{\omega}\right) \left(\displaystyle \sum_{i=0}^k\; {{\rm Ric}(\omega_{\varphi})}^{i}\wedge\omega^{k-i} \right)
    \wedge {\omega_{\varphi}}^{n-k}.
\]
\end{defi}
If there is no possible confusion, we will often drop the subscript $\omega$ in the following.

\begin{rem} If $k = n=1,$ then
\[
  E_1^0 = {1\over V}\; \displaystyle \int_M\;  \left( \log {{\omega_{\varphi}} \over \omega}- h_{\omega}\right)
  \;\left(R \; {\omega_{\varphi}} + 1\cdot \omega\right).
\]

This is analogous to the well-known Liouville energy on Riemann
surfaces.
\end{rem}

Next for each $k=0,1,2,\cdots, n-1$, we will define $J_{k, \omega}$ as follows:
Let $\varphi(t) $ ($t\in [0,1]$) be a path from $0$ to $\varphi$ in
${\cal P}(M,\omega)$, we define
\[ J_{k,\o}(\varphi) = -{n-k\over V}
\int_0^1 \int_M {{\partial \varphi}\over{\partial t}} \left({\omega_{\varphi}}^{k+1} -
{\omega}^{k+1}\right)\wedge {\omega_{\varphi}}^{n-k-1}\wedge dt.
\]

One can show that the integral on the right is independent of
choices of the path. This is because ${\cal P}(M,\omega)$ is
simply-connected and its derivative has nothing to do with the
path. Clearly, we have
\[
{{d J_{k,\omega}}\over {d\,t}} = - {{n-k}\over V}\;\displaystyle \int_M\;
{{\partial \varphi}\over{\partial t}} \; \left({\omega_{\varphi}}^{k+1} -
{\omega}^{k+1}\right)\wedge {\omega_{\varphi}}^{n-k-1},
\]

For simplicity, we will often drop the subscript $\o$ in the following.

\begin{rem} If $k= n-1$, then
\[
  {{d J_{n-1}}\over {d\,t}} = - {1\over V}\;\displaystyle \int_M\; {{\partial \varphi}\over{\partial t}}
  \; \left({\omega_{\varphi}}^{n} - {\omega}^{n}\right)
\]

Thus  $J_{n-1,\omega} = J_{\omega}$ is just the generalized
energy functional (cf. \cite{DingTian93} and also Section 3.1).
\end{rem}

\begin{prop} For each $k=0, 1, \cdots, n-1$, we have the following explicit formula for $J_k$:
\begin{equation}
J_k(\varphi) = {{n-k}\over V}\;\displaystyle \sum_{j=0}^{n-k-1}
\displaystyle \sum_{i=0}^k
 \displaystyle \sum_{s=0}^{n-i-j-1}  \;c_{sij} \;\displaystyle \int_M\; \sqrt{-1} \;
\partial \varphi \wedge \overline{\partial}\varphi  \wedge
{\omega_{\varphi}}^s \wedge \omega^{n-1-s} , \label{eq:jk}
\end{equation}
where $c_{sij} $ is
\[
{{(-1)^{n-i-j-s-1} }\over {(n-i-j+1)}}\; {{k+1} \choose i}
{{n-k-1} \choose {j}} {{n-i-j-1}\choose s}
\]
\end{prop}
\begin{proof} We will calculate $J_k (\varphi) $ via a trivial path $
t\varphi \in {\cal P}(M,\omega)$ (the corresponding K\"ahler
metrics are $\omega + t \sqrt{-1} \partial \overline{\partial}
\varphi$).

\[
\begin{array}{ll} & \qquad J_k(\varphi) \\
  & =  {k-n\over V}\displaystyle \int_0^1 \displaystyle
 \int_M\varphi\left(
\omega_{t\varphi}^{k+1} - {\omega}^{k+1}\right)
\wedge \omega_{t\varphi}^{n-k-1} \wedge dt\\
&  =  {{k-n}\over V}\displaystyle \int_0^1 \displaystyle \int_M
\varphi \displaystyle \sum_{i=0}^k  {{k+1} \choose i}
 \omega^i \wedge (\sqrt{-1}\partial
\overline{\partial} \varphi )^{k+1-i} t^{k+1-i}
 \\
 & \qquad \qquad \wedge \displaystyle \sum_{j=0}^{n-k-1} {{n-k-1}
\choose j}
 \omega^j\wedge (\sqrt{-1}\partial \overline{\partial} \varphi )^{n-k-1-j} t^{n-k-1-j} \wedge d\,t\\
 & = {{k-n}\over V}\displaystyle \int_0^1 \displaystyle \int_M
\varphi \displaystyle \sum_{i=0}^k  \displaystyle
\sum_{j=0}^{n-k-1} {{k+1} \choose i} {{n-k-1} \choose j}
 \omega^{i+j} \wedge (t \sqrt{-1}\partial \overline{\partial} \varphi \}^{n-i-j} \wedge dt\\
 & =  -{n-k\over V}\displaystyle \int_M \varphi \; \displaystyle
\sum_{j=0}^{n-k-1} \displaystyle \sum_{i=0}^k {1 \over {n-i-j+1}}
{{k+1} \choose i} {{n-k-1}\choose {j}} \omega^{i+j}\wedge
 (\sqrt{-1}\partial \overline{\partial} \varphi )^{n-i-j}\\
 & =   \displaystyle \sum_{i=0}^k \displaystyle \sum_{j=0}^{n-k-1} {{k+1} \choose i}
{{n-k-1} \choose {j}} {{(n-k)\sqrt{-1}}\over {(n-i-j+1)V}}
\displaystyle \int_M\;\partial \varphi \wedge
\overline{\partial}\varphi \wedge \omega^{i+j}\wedge
(\omega_{\varphi} - \omega )^{n-i-j-1}\\
 & = \displaystyle
\sum_{j=0}^{n-k-1} \displaystyle \sum_{i=0}^k \displaystyle
\sum_{s=0}^{n-i-j-1} \; {{(n-k)(-1)^{n-i-j-s-1} }\over
{(n-i-j+1)V}}\; {{k+1} \choose i} {{n-k-1} \choose {j}}
{{n-i-j-1}\choose s}
\\  & \qquad \qquad \qquad \qquad \qquad \qquad \displaystyle \int_M
\; \sqrt{-1} \;\partial \varphi \wedge \overline{\partial}\varphi
\wedge \omega^{i+j} \wedge {\omega_{\varphi}}^s\wedge
\omega^{n-s-i-j-1}.
\end{array}
\]
\end{proof}

The following is an immediate corollary of Formula (\ref{eq:jk}).
\begin{cor}
For each $k$, there is a uniform constant $a_k$ such that for any $\varphi \in {\cal P} (M,\omega)$,
\[
  | J_{k,\omega}(\varphi) | \leq a_k \cdot J_{\omega}(\varphi).
\]
\end{cor}

This follows from Formula (\ref{eq:jk}) and the explicit
expression: \[ J_{\omega}(\varphi) = {1 \over V} \displaystyle
\sum_{i=0}^{n-1}\;\displaystyle\int_M\;{{i+1}\over{n+1}}\;\sqrt{-1}
\;
\partial \varphi \wedge \overline{\partial}\varphi \wedge \omega^i \wedge {\omega_{\varphi}}^{n-1-i}.
\]

Now we simply define $E_{k,\omega}$ ($k=0,1,\cdots, n$) by
$$E_{k,\omega}(\varphi) = E_{k,\omega}^0 (\varphi) - J_{k,\omega}(\varphi),$$
where we set $J_{n,\omega}=0$.

In the following, we will often write $E_k$ for $E_{k,\omega}$ if no confusion may occur.

\subsection{Derivative of $E_k$}

In this subsection, we derive a few basic properties of $E_k$.

\begin{theo} For any $k=0,1,\cdots, n$, we have\footnote{In a non
canonical K\"ahler class, we need to modify the definition
slightly since $h_{\omega}$ is not defined there. For any
$k=0,1,\cdots, n,\;$ we define
\[
\begin{array}{lcl}
  E_{k,\omega}(\varphi) & = &
 {1\over V}\; \displaystyle \int_M\;   \log {{\omega_{\varphi}}^n \over \omega^n}
\; \left(\displaystyle \sum_{i=0}^k\; {{\rm
Ric}(\omega_{\varphi})}^{i}\wedge {\rm Ric}(\omega)^{k-i} \right)
    \wedge {\omega_{\varphi}}^{n-k}  \\
    &  & \qquad \qquad  -
    {{n-k}\over V} \displaystyle \int_M\; \varphi
    \left({\rm Ric}(\omega)^{k+1} - \omega^{k+1}\right) \wedge \omega^{n-k-1}
    - J_{k,\omega}(\varphi).
    \end{array}
\]
The second integral on the right is to offset the change from
$\omega$ to $Ric(\omega)$ in the first term. The derivative of
this functional is exactly same as in the canonical K\"ahler
class. In other words, the Euler-Lagrange equation is not
changed.}
\begin{eqnarray}
{d E_k \over dt} & = & {{k+1}\over V} \displaystyle \int_M
\Delta_{\varphi}\left( {{\partial \varphi}\over {\partial t}}
\right )\; {{\rm Ric}(\omega_{\varphi})}^{k} \wedge
{\omega_{\varphi}}^{n-k}
\nonumber \\
& & \qquad -
 {{n-k}\over V}\displaystyle \int_M {{\partial \varphi}\over {\partial t}} \left({{\rm Ric}(\omega_{\varphi})}^{k+1}
 - {\omega_{\varphi}}^{k+1}\right) \wedge  {\omega_{\varphi}}^{n-k-1}.
\label{eq:decay functional0}
\end{eqnarray}
Here $\{\varphi(t)\}$ is any path in ${\cal P}(M,\omega)$.
\end{theo}
\begin{rem} When $k=0,$ we have
\[
{d E_0 \over dt} ={{n}\over V}\displaystyle \int_M {{\partial
\varphi}\over {\partial t}} \left({\rm Ric}(\omega_{\varphi})
 - \omega_{\varphi}\right) \wedge  {\omega_{\varphi}}^{n-1}.
\]
Thus $E_0$ is a multiple of the well known K-energy function
introduced by T. Mabuchi.
 \end{rem}
\begin{proof} We suffice to compute the derivatives of $E_k^0$ ($k=0,1,\cdots, n$).
Put $F =  \log {\omega_{\varphi}^n \over {\omega}^n } - h_{\omega}$. It is clear that
\[
\begin{array}{lcl}
\sqrt{-1} \partial \overline{\partial} F & = & {\rm Ric}(\omega) - {\rm Ric}(\omega_{\varphi}) -
\sqrt{-1} \partial \overline{\partial} h_{\omega}\\
  &  = &\omega-{\rm Ric}(\omega_{\varphi})
\end{array}
\]
and
\[
  {{\partial {\rm Ric}(\omega_{\varphi})} \over {\partial t}} =
- \sqrt{-1} \partial \overline{\partial} \Delta_{\varphi}\left ( {{\partial \varphi}\over {\partial t}}\right ).
\]

Thus,
\[
\begin{array}{lcl}
& & {{d E_k^0}\over{d t}} \\
& = & {1\over V}\displaystyle \int_M
\Delta_{\varphi} \left ({{\partial \varphi}\over {\partial
t}}\right ) \left(\displaystyle \sum_{i=0}^k {{\rm
Ric}(\omega_{\varphi})}^{i}\wedge\omega^{k-i} \right) \wedge
{\omega_{\varphi}}^{n-k}
\\ & & +  {{1}\over V} \displaystyle \int_M F \displaystyle \sum_{i=0}^k i\; {{\rm Ric}(\omega_{\varphi})}^{i-1}
\wedge\omega^{k-i} \wedge \left( - \sqrt{-1} \partial
\overline{\partial} \Delta_{\varphi}
 {{\partial \varphi}\over {\partial t}}\right) \wedge {\omega_{\varphi}}^{n-k} \\
& & \qquad +
 {{n-k}\over V} \displaystyle \int_M  F \left(\displaystyle \sum_{i=0}^k
{{\rm Ric}(\omega_{\varphi})}^{i}\wedge\omega^{k-i} \right) \wedge
\left( \sqrt{-1} \partial \overline{\partial} \left({{\partial
\varphi}\over {\partial t}}\right)\right)\wedge
{\omega_{\varphi}}^{n-k-1} \\
& = & {1\over V} \displaystyle \int_M \Delta_{\varphi}\left( {{\partial \varphi}\over {\partial t}}\right )
 \left (\displaystyle \sum_{i=0}^k {{\rm Ric}(\omega_{\varphi})}^{i}\wedge\omega^{k-i} \wedge
{\omega_{\varphi}}^{n-k}\right. \\ & & - \left. \sqrt{-1}
\partial \overline{\partial} F \wedge \displaystyle \sum_{i=0}^k
\;i\;
{{\rm Ric}(\omega_{\varphi})}^{i-1}\wedge\omega^{k-i}\wedge {\omega_{\varphi}}^{n-k} \right) \\
& & \qquad +
 {{n-k}\over V}\displaystyle \int_M {{\partial \varphi}\over {\partial t}}
\left(\displaystyle \sum_{i=0}^k {{\rm Ric}(\omega_{\varphi})}^{i}\wedge\omega^{k-i} \right) \wedge
\left( \sqrt{-1} \partial \overline{\partial} F \right)\wedge  {\omega_{\varphi}}^{n-k-1}.
\end{array}
\]

Plugging $\sqrt{-1} \partial \overline{\partial} F = \omega -
{\rm Ric}(\omega_{\varphi})$, we obtain
\[
\begin{array}{lcl}
& & {{d E_k^0}\over{d t}} \\
& = & {1\over V} \displaystyle \int_M
\Delta_{\varphi}\left( {{\partial \varphi}\over {\partial
t}}\right ) \left (\displaystyle \sum_{i=0}^k {{\rm
Ric}(\omega_{\varphi})}^{i}\wedge\omega^{k-i} \wedge
{\omega_{\varphi}}^{n-k} \right. \\
& & \qquad + \left.\left({\rm Ric}(\omega_{\varphi}) -\omega
\right)\wedge \displaystyle \sum_{i=0}^k\; i \;{{\rm
Ric}(\omega_{\varphi})}^{i-1}\wedge\omega^{k-i} \wedge
{\omega_{\varphi}}^{n-k}\right)
\\ & & \qquad + {{n-k}\over V} \displaystyle \int_M
{{\partial \varphi}\over {\partial t}} \left(\displaystyle \sum_{i=0}^k {{\rm Ric}(\omega_{\varphi})}^{i}\wedge
\omega^{k-i} \right) \wedge \left(
 \omega -{\rm Ric}(\omega_{\varphi}) \right)\wedge  {\omega_{\varphi}}^{n-k-1}.
\end{array}
\]

Now we recall a polynomial identity: For any two variables $x,y$
and any integer $k>0$, we have
\begin{equation}
\displaystyle \sum_{i=0}^k \; x^i \;y^{k-i} + (x-y) \displaystyle
\sum_{i=0}^k \;i x^{i-1} y^{k-i} =  (k+1) \;x^{k}.
\label{eq:polynomial}
\end{equation}

Applying this identity to the first integral above, we get
\[
\begin{array}{lcl}
{{d E_k^0}\over{d t}}
& = &{k+1\over V}\displaystyle \int_M \Delta_{\varphi} \left( {{\partial \varphi}\over {\partial t}}\right )
{{\rm Ric}(\omega_{\varphi})}^{k}\wedge{\omega_{\varphi}}^{n-k} \\
& & \qquad +
 {{n-k}\over V} \displaystyle \int_M {{\partial \varphi}\over {\partial t}}
\left(\omega^{k+1} - {{\rm {\rm Ric}}(\omega_{\varphi})}^{k+1}\right)\wedge {\omega_{\varphi}}^{n-k-1}.
\end{array}
\]

The theorem follows from this and explicit expression of the
derivative of $J_k$.
\end{proof}

From this theorem, we can show that all $E_k$ satisfy a cocycle condition.

\begin{cor}
For each $k=0,1,\cdots,n$, the functional $E_{k,\omega}$ satisfies the following:
For any $\varphi$ and $\psi$ in ${\cal P}(M,\omega)$,
\[
E_{k,\omega} (\varphi) +  E_{k,\omega_\varphi}(\psi -\varphi) =
E_{k,\omega}(\psi).
\]
\end{cor}

Let us write down the Euler-Lagrange equation for the functional
$E_k$($k=0,1,\cdots, n$). Recall the expansion formula
(\ref{symm}) in $t$:
$$ \left ( \omega_\varphi + t \;{\rm Ric} (\omega_\varphi) \right
)^n = \left( \sum_{k=0}^n \sigma_k(\omega_\varphi) t^k \right )
\omega_\varphi^n.$$ Clearly, $\sigma_0(\omega_\varphi)=1$,
$\sigma_1(\omega_\varphi) = R(\omega_\varphi)$, the scalar
curvature of $\omega_\varphi$. In general, $\sigma_k$ is a k-th
symmetric polynomial of Ricci curvature. The Euler-Lagrange
equation of $E_k$ is
\[
(k+1)  \Delta_\varphi \sigma_k (\omega_\varphi)
-(n-k)\sigma_{k+1}(\omega_\varphi ) = c_k,
\]
where $\Delta_\varphi$ is the Laplacian of the metric
$\omega_\varphi$ and $c_k$ is the constant $$ -(n-k)
c_1(M)^{k+1}\cup [\omega]^{n-k-1} ([M]).$$ Clearly,
K\"ahler-Einstein metrics are solutions to the above equation for
any $k$. If the K\"ahler class is canonical, one can show that
for $k=n$, K\"ahler-Einstein metrics are the only solutions of the
Euler-Lagrange equation with positive Ricci curvature. However,
It is not clear what the critical points are in other K\"ahler
classes. But it certainly merit further study of these equations.

\begin{prop} Along the K\"ahler Ricci flow, we have
\begin{equation}
{d E_k \over d t} \leq - {{k+1}\over V} \displaystyle \int_M
(R(\omega_\varphi)-r) {\rm Ric}(\omega_{\varphi})^{k} \wedge
{\omega_{\varphi}}^{n-k}. \label{eq: Ekdecreases}
\end{equation}
When $k=0 , 1$, we have
\begin{eqnarray}
{{d E_0 }\over{d\,t}} & =  & -{{n\sqrt{-1}}\over V} \displaystyle
\int_M \partial {{\partial \varphi}\over {\partial t}} \wedge
\overline{\partial} {{\partial \varphi}\over{\partial
t}}  {\omega_{\varphi}}^{n-1}\leq 0,\\
{{d E_1 }\over{d t}} & \leq & - {{2}\over V} \displaystyle \int_M
(R(\omega_\varphi)-r)^2 {\omega_{\varphi}}^{n} \leq  0. \nonumber
\end{eqnarray}
In particular, both $E_0$ and $E_1$ are decreasing along the
K\"ahler Ricci flow.
\end{prop}
\begin{proof} Along the K\"ahler Ricci flow, we have
\[\Delta_{\varphi} \left ({{\partial \varphi}\over {\partial t}} \right ) = r - R(\omega_\varphi).
\]
Here $r$ is again the average of the scalar curvature
$R(\omega_\varphi)$. We also have
\[
\begin{array}{lcl}
  \sqrt{-1} \partial \overline{\partial} {{\partial \varphi}\over{\partial t}} &  = &
  \sqrt{-1} \partial \overline{\partial} \left(\log {{\omega_{\varphi}}^n \over {\omega}^n } + \varphi - h_{\omega}
 \right) \\
& = & - {\rm Ric}(\omega_{\varphi}) + \left(\omega - \sqrt{-1}
\partial \overline{\partial}
h_{\omega} \right) + \sqrt{-1} \partial \overline{\partial} \;\varphi\\
 & = & - {\rm Ric}_(\omega_{\varphi}) + \omega + \sqrt{-1} \partial \overline{\partial} \;\varphi =
 \omega_{\varphi} -{\rm Ric}(\omega_{\varphi}).
\end{array}
\]
Therefore,
\begin{eqnarray}
& & {{d E_k }\over{d t}} \\
& =  & - {{k+1}\over V} \displaystyle
\int_M (R(\omega_\varphi)-r)
 {{\rm Ric}(\omega_{\varphi})}^{k} \wedge {\omega_{\varphi}}^{n-k}
\nonumber \\ &  &  + {{n-k}\over V} \displaystyle
\int_M\;{{\partial \varphi}\over {\partial t}}
 \sqrt{-1} \partial \overline{\partial}
 \left({{\partial \varphi}\over{\partial t}}\right)
 \wedge \displaystyle \sum_{i=0}^{k+1} {\rm Ric}(\omega_{\varphi})^{k+1-i}\wedge  \omega_{\varphi}^{n-k+i-1} \\
&\leq & - {{k+1}\over V} \displaystyle \int_M
(R(\omega_\varphi)-r) {{\rm Ric}(\omega_{\varphi})}^{k} \wedge
{\omega_{\varphi}}^{n-k}.
\end{eqnarray}
\end{proof}

The following is an easy corollary of the above, but it will be
crucial in our proof.

\begin{theo} Let $\varphi(t)$ be the global solution of the K\"ahler Ricci flow.
Then for any $T>0$, we have
\[
{{k+1}\over V}\displaystyle \int_0^T \displaystyle \int_M\
(R(\omega_{\varphi})-r)\; {{\rm Ric}(\omega_{\varphi})}^{k}
\wedge {\omega_{\varphi}}^{n-k} \; d\,t \le E_k(\varphi(0))  -
E_k(\varphi(T)).
\]
When $k=1,$ this reduces to
 \[
{{2}\over V}\displaystyle \int_0^T \displaystyle \int_M\
(R(\omega_{\varphi})-r)^2 {\omega_{\varphi}}^{n} \; d\,t \le
E_1(\varphi(0)) - E_1(\varphi(T))
 \]
 In particular, if
$E_k(\varphi(t))$ is uniformly bounded from below, then for any
sequence of positive numbers $\epsilon_i$ with $\displaystyle
\lim_{i\to \infty} \epsilon_i =0$, there exists a sequence of
$t_i$ such that
\[
\displaystyle \sum_{k=0}^n {{k+1}\over V} \displaystyle \int_M
(R(\omega_{\varphi(t_i)})-r)
 {{\rm Ric}(\omega_{\varphi(t_i)})}^{k} \wedge {\omega_{\varphi(t_i)}}^{n-k}
\le \epsilon_i.
\]

When $k=1$, this becomes
\[
{{1}\over V} \displaystyle \int_M (R(\omega_{\varphi(t_i)})-r)
^2\; {\omega_{\varphi(t_i)}}^{n} \le \epsilon_i.
\]
\end{theo}

In order to have integral bounds of curvature from these
inequalities, we need to bound these functionals $E_k$ from below.
The following provides a way of achieving it.

\begin{lem}
Let $\varphi$ be in ${\cal P}(M,\omega)$ such that ${\rm
Ric}(\omega_\varphi) \ge 0$. Then there is a uniform constant
$c=c(\omega)$ such that
\[
E_k(\varphi) \ge - e^{c\left(1 + {\rm max} \{0,
\nu_\omega(\varphi)\} + J_\omega(\varphi)\right)}.
\]
\end{lem}

\begin{proof} We will always denote by $c$ a constant depending
only on $\omega$. By the definition of $E_k$ and Corollary 4.5,
we have
\[
E_k \ge  {1\over V} \displaystyle \int_M \left( \log
{\omega_{\varphi}^n \over \omega^n} \right) \left(\displaystyle
\sum_{i=0}^k\; {{\rm Ric}(\omega_{\varphi})}^{i}\wedge\omega^{n-i}
\right) \,-\, c \left (1+J_\omega(\varphi)\right ).
\]
In particular, since $E_0$ is just the K-energy, we have
\[
{1\over V} \displaystyle \int_M \left( \log {\omega_{\varphi}^n
\over \omega^n} \right) \omega^n_\varphi \le \nu_\omega(\varphi)
+ c (1+J_\omega(\varphi)).
\]
Then the lemma follows from the above two inequalities and the
volume estimate in Proposition 3.13.
\end{proof}

Because of the monotonicity of the K-energy along the K\"ahler
Ricci flow, the K-energy $\nu_\omega(\varphi)$ is bounded. Hence,
in order to bound $E_k$, we suffice to bound the generalized
energy $J_\omega(\varphi)$ along the K\"ahler Ricci flow. The
trouble is that $J_\omega(\varphi)$ may not be bounded along the
flow. We will bound $J_\omega$ for modified K\"ahler Ricci flow,
which turns out to be sufficient (cf. Section 6).

\section{New holomorphic invariants}

In this section, we want to show that on any K\"ahler-Einstein
manifolds, $E_k$ ($k=0,1,\cdots, n$) are invariant under
automorphisms. First we want to show that the derivatives of
$E_k$ in the direction of holomorphic vector field give us
holomorphic invariants of the K\"ahler class.

Let $X$ be a holomorphic vector field and $\omega$ be a K\"ahler
metric. Then $i_X\omega$ is a $\overline \partial $-closed (0,1)
form, by the Hodge theorem, we can decompose $i_X\omega$ into a
parallel $\alpha_X$ form plus $\sqrt{-1} \overline \partial
\theta_X$, where $\theta_X$ is some function. For simplicity, we
will assume that $\alpha_X=0$. This is automatically true if $M$
is simply-connected. We will call that $\theta_X$ is a potential
of $X$ with respect to $\omega$. It is unique modulo addition of
constants. Note that $L_X (\omega) = \sqrt{-1} \partial
\overline{\partial} \theta_X$. Now we define $\Im_k (X,\omega)$
for each $k=0,1,\cdots, n$ by
\[
\begin{array}{ll}
\Im_k (X,\omega)  = (n-k) \displaystyle \int_M\theta_X\;
{\omega}^n\\
\qquad  + \displaystyle \int_M \left( (k+1) \Delta \theta_X
\;{{\rm Ric}(\omega)}^{k}\wedge {\omega}^{n-k} - (n-k)\; \theta_X
\;{{\rm Ric}(\omega)}^{k+1} \wedge {\omega}^{n-k-1}\right ).
\end{array}
\]
Here and in the following, $\Delta$ denotes the Laplacian of
$\omega$. Clearly, the identity is unchanged if we replace
$\theta_X$ by $\theta_X + c$ for any constant $c$.

The next theorem assures that the above integral gives rise to a
holomorphic invariant.

\begin{theo} The integral $\Im_k (X,\omega)$ is independent of choices of
K\"ahler metrics in the K\"ahler class $[\omega]$, that is,
$\Im_k (X,\omega)=\Im_k (X,\omega')$ so long as the K\"ahler
forms $\omega$ and $\omega'$ represent the same K\"ahler class.
Hence, the integral $\Im_k (X,\omega)$ is a holomorphic
invariant, which will be denoted by $\Im_k (X,[\omega])$.
\end{theo}

\begin{rem}
When $k=0$, we have
\[\begin{array}{lcl}
\Im_0 (X,\omega)& = & \displaystyle \int_M \Delta \theta_X \;
{\omega}^{n} + n \,\theta_X
\left( {\omega}- {{\rm Ric}(\omega)}\right) \wedge {\omega}^{n-1}\\
& = & - n \displaystyle \int_M \theta_X
\;\Delta\;h_{\omega}\;{\omega}^{n} =  n \displaystyle \int_M
X(h_{\omega})\;{\omega}^{n}.
\end{array}
\]
Thus $\Im_0 (X,[\omega])$ is a multiple  of the Futaki invariant
(\cite{futaki83}).
\end{rem}

If $[\omega]$ is a canonical K\"ahler class and there is a
K\"ahler-Einstein metric on $M$, then we can choose $\omega$ such
that ${\rm Ric}(\omega) = \omega$ and deduce
\[
\Im_k (X,\omega) = (k+1) \displaystyle \int_M  \Delta \theta_X
\;\omega^n = 0.
\]
Therefore, we have

\begin{cor}
The above invariants $\Im_k (X, c_1(M))$ all vanish for any
holomorphic vector fields $X$ on a compact K\"ahler-Einstein
manifold. In particular, these invariants all vanish on $\CC P^n$.
\end{cor}

Before we prove this theorem, we first use to show the invariance
of $E_k$ under automorphisms.

\begin{prop} Let $X$ be a holomorphic vector field and $\{\Phi(t)\}_{|t| <
\infty}$ be the one-parameter subgroup of automorphisms induced by
$Re(X)$. Then
\[
{d E_k (\varphi_t)\over dt}  = { 1 \over V} {\rm Re}(\Im_{k}(X,
\omega)),
\]
where $\varphi_t$ are the K\"ahler potentials of
$\Phi_t^*\omega$, i.e., $\Phi_t^*\omega=\omega + \sqrt{-1}
\partial\overline\partial \varphi_t$.
\end{prop}
\begin{proof} Differentiating $\Phi_t^*\omega = \omega + \sqrt{-1}
\partial\overline\partial \varphi_t$, we get
\[
L_{Re(X)}\omega = \sqrt{-1}
\partial\overline\partial \left ( {\partial \varphi_t \over \partial t}\right ).
\]
On the other hand, since $L_X\omega = \sqrt{-1}
\partial\overline\partial \theta_X$, we have
\[
{\partial \varphi_t \over \partial t} = Re(\theta_X) + c,
\]
where $c$ is a constant. It follows
 \[
\begin{array}{lcl}
{{d E_k }\over{d t}}&=& {1\over V} \displaystyle \int_M \left (
(k+1) \Delta \left({\partial \varphi_t \over \partial t}\right
)\,{{\rm Ric}(\omega)}^{k} \wedge {\omega}^{n-k} \right.\\
&&\qquad \left. -(n-k) \,{\partial \varphi_t \over \partial t}
\left({{\rm Ric}(\omega)}^{k+1} - {\omega}^{k+1}\right) \wedge
{\omega}^{n-k-1}
\right )  \\
& = &  {1\over V}\; Re( \Im_{k}(X,\omega))\\
& = & 0.
\end{array}
\]
\end{proof}
An immediate corollary is

\begin{cor} On a K\"ahler-Einstein manifold $M$ with $c_1(M)=[\omega]$,
all functionals $E_{k,\omega}$ ($k=0,1,\cdots,n$) are invariant
under automorphisms of $M$.
\end{cor}

\begin{rem}
It also follows from the above proposition that $E_k$ has a lower
bound only if $\Im_{k}(X, \omega) = 0$.
\end{rem}

The rest of this section is devoted to proving this theorem. We
will follow the arguments in \cite{tian98}. For this purpose, we
first formulate $\Im_k$ in terms of some particular forms, i.e.,
the Bott-Chern forms.

\begin{lem}
There exists a matrix $\left(c_{ij}\right)$ ($1\le i, j\le n+1$)
such that
\[
\begin{array}{lcl}
\Im_{k-1}(X,\omega) &=& - {{n-k+1 +\upsilon_k}\over {n+1}}
\displaystyle \int_M \left(-\theta_X +  \omega\right)^{n+1} \\
& & \qquad + {1 \over {n+1 \choose k}} \displaystyle
\sum_{i=1}^{n+1} c_{i k} \displaystyle \int_M \left(-\theta_X +
\omega + i (\Delta \theta_X +
   {\rm Ric}(\omega)) \right)^{n+1},
\end{array}
\]
where the matrix $\left(c_{ij}\right)$ is the inverse matrix of
the well known Vandermonde matrix \footnote{There is an explicit
way of finding $c_{ij}$, which we learned from E. Calabi. Let us
define a sequence of polynomials of degree $n+1$ by
\[
f_i(x) = \displaystyle \sum_{j=1}^{n+1} \; c_{i j} x^j,\qquad
\forall\; i=1,2,\cdots, n+1.
\]
Since $(c_{ij}) $ is the inverse matrix of Vandermonde matrix:
\[
\left(\begin{array} {llll}   1 & 2 & \cdots & n+1 \\
                               1^2 & 2^2 & \cdots & (n+1)^2 \\
                               1^3 & 2^3 & \cdots & (n+1)^3 \\
                               \vdots & \vdots & \ddots &\vdots \\
                              1^{n+1} & 2^{n+1} & \cdots & (n+1)^{n+1}\end{array}
\right),
\]
we obtain
\[
  f_i(k) = \displaystyle \sum_{j=1}^{n+1}  c_{i j} k^j = \delta_{k j},\qquad \forall\; i,k=1,2,\cdots, n+1.
\]
It follows that for each $i=1,2,\cdots, n+1$,
\[
  f_i(x) =  \displaystyle \sum_{j=1}^{n+1} \; c_{i j} x^j = {{x\;\Pi_{k\neq i, k = 1}^{n+1} (x-k) }
  \over {i\;\displaystyle \Pi_{k\neq i, k = 1}^{n+1} (i-k) }}.
\]
Thus $c_{ij}$ can be found.} and
\[ \upsilon_k = {{n+1}\over{n+1\choose k}} \displaystyle\sum_{i=1}^{n+1}
\;c_{ik}.\qquad where \qquad k=1,2,\cdots, n+1.
\]
\end{lem}

\begin{proof} Consider
\[
\begin{array}{lcl}
I_{pq} & = & \displaystyle \int_M\; \left(- p\; \theta_X + q\;
\Delta \theta_X \right)\;
  \left( q \;{\rm Ric}(\omega) + p\; \omega \right)^n \\
& = & {1\over {n+1}} \displaystyle\; \int_M\; \left(- p \;\theta_X
+ q \;\Delta \theta_X+
  q \;{\rm Ric}(\omega) + p\; \omega \right)^{n+1}.
\end{array}
\]
Note that the only forms of degree $2n$ contribute to the above
integral.

Expanding the integrand, we have
\[
\begin{array}{lcl}
& & I_{pq} \\
& = & \displaystyle \int_M \left(- p\; \theta_X + q\;
\Delta \theta_X \right) \left( \displaystyle\;\sum_{k=0}^n q^k
p^{n-k}{n\choose k}\; {{\rm Ric}(\omega)}^k\wedge {\omega}^{n-k}
 \right)\\
& = & - \displaystyle \int_M \theta_X\left( \displaystyle
\sum_{k=0}^n q^k p^{n-k + 1 }
{n\choose k}\;{{\rm Ric}(\omega)}^k\wedge {\omega}^{n-k} \right) \\
 & & \qquad \qquad + \displaystyle \int_M \Delta \theta_X\left( \displaystyle\sum_{k=0}^n q^{k+1}p^{n-k}
 {n\choose k}\;{{\rm Ric}(\omega)}^k \wedge{\omega}^{n-k} \right) \\
& = & -  \displaystyle \int_M \theta_X\left(
\displaystyle\sum_{k=0}^n q^k p^{n-k + 1 }
{n\choose k}\;{{\rm Ric}(\omega)}^k \wedge {\omega}^{n-k} \right) \\
 & & \qquad \qquad + \displaystyle \int_M \Delta\theta_X\left( \displaystyle\sum_{k=1}^{n+1} q^{k}p^{n-k+1}
 {n\choose k-1}\;{{\rm Ric}(\omega)}^{k-1}\wedge {\omega}^{n-k+1} \right) \\
& = & \displaystyle \sum_{k=0}^{n+1} q^{k}
p^{n-k+1} {{n!} \over {k! (n-k+1)!}} \\
& & \; \cdot\displaystyle\int_M \left( k\; \Delta\theta_X\;{{\rm
Ric}(\omega)}^{k-1} \wedge {\omega}^{n-k+1} - (n-k+1)\;\theta_X
\;{{\rm Ric}(\omega)}^k \wedge{\omega}^{n-k}\right).
\end{array}
\]

Now set $p=1$ and observe ($k=0,1,2,\cdots, n$)
\[
\begin{array}{lcl}
&&\Im_k(X,\omega) - (n-k) \displaystyle \int_M \theta_X
{\omega}^{n} \\
&=& \displaystyle \int_M\left((k+1) \Delta\theta_X\;{{\rm
Ric}(\omega)}^{k}
 \wedge {\omega}^{n-k} - (n- k )\;\theta_X {{\rm Ric}(\omega)}^{k+1}\wedge {\omega}^{n-k-1}\right).
 \end{array}
\]

Then
\[
\begin{array}{lcl}
I_{1q} &  = &  -\displaystyle \int_M \theta_X\;\omega^n \\
& & \qquad + {1\over (n+1)}\displaystyle\sum_{k=1}^{n+1}
{n+1\choose k}\; \left({\cal F}_{k-1}(X,\omega)- (n-k+1)
\displaystyle \int_M \theta_X {\omega}^{n}\right)\;q^{k} ,
\end{array}
\]
or equivalently,
\[
{1\over (n+1)}\displaystyle\sum_{k=1}^{n+1} {n+1\choose
k}\;\left( \Im_{k-1}(X,\omega) - (n-k+1) \displaystyle \int_M
\theta_X {\omega}^{n}\right) \;q^{k} = I_{1q} + \displaystyle
\int_M \theta_X {\omega}^{n}.
\]

Since $\left(c_{ij}\right)$ is the inverse matrix of the
Vandermonde matrix, we have
\[
\begin{array}{lcl}
&&{ {n+1\choose k}\over (n+1)} \; \left(\Im_{k-1}(X,\omega) -
(n-k+1) \displaystyle \int_M \theta_X {\omega}^{n}\right)\\
&=&\displaystyle \sum_{i=1}^{n+1} c_{ik} \left(I_{1 i}+
\displaystyle
\int_M \theta_X {\omega}^{n}\right)\\
& = & \displaystyle \sum_{i=1}^{n+1} c_{ik} I_{1i} +
\displaystyle \sum_{i=1}^{n+1} c_{ik}\;\displaystyle \int_M
\theta_X {\omega}^n \\
& = & \displaystyle \sum_{i=1}^{n+1} c_{ik} I_{1i} + \upsilon_k {
{n+1\choose k}\over (n+1)}\displaystyle \int_M \theta_X
{\omega}^n.
\end{array}
\]

The lemma follows from this since
\[- \displaystyle \int_M
\theta_X {\omega}^n = {1\over n+1} \displaystyle \int_M \left (
\theta_X + {\omega}\right )^{n+1}. \]
\end{proof}

Now we continue the proof of Theorem 5.1. We suffice to prove the
independence of $I_{pq}$. First we observe that
\[
\sqrt{-1}\;\overline{\partial} \theta_X = i_X \omega~~~{\rm
and}~~~ \sqrt{-1}\;\overline{\partial} \Delta \theta_X = - i_X
{\rm Ric}(\omega).
\]
The second identity can be checked as follows: Suppose $\omega
=\sqrt{-1} \displaystyle \sum_{i,j=1}^n\; g_{i \overline{j}}\;
dz_i \wedge d \overline{z_j} $ in local coordinates. Then
\[
\begin{array}{lcl} i_X {\rm Ric}(\omega) &= & -\sqrt{-1}\; \b \p \Bl X^i {\p \over {\p z_i}}
                     \log \det (g_{k\b l})\Br\\
                  &= &-\sqrt{-1}\;\b \p \Bl X^i g^{k\b l} {{\p g_{i \b l}} \over {\p z_k}}
                     \Br \\
                  &= & -\sqrt{-1}\;\b \p \Bl g^{k \b l} {\p \over {\p z_k}} (X^i g_{i \b l})-
                    g^{k \b l} g_{i \b l} {{\p X^i} \over {\p z_k}}\Br \\
                  &= & -\sqrt{-1}\;\b \p \Bl g^{k \b l} {\p \over {\p z_k}} (X^i g_{i\b l})
                     \Br \\
                  &= & -\sqrt{-1} \;\b \p \Delta _g \th _X.
\end{array}
\]

Since the space of K\"ahler metrics is path-connected, it suffices
to show that $I_{pq}$ is invariant when we deform the K\"ahler
potential along any path $\varphi_t\in {\cal P}(M,\omega)$. To
emphasis the dependence on $\omega_\varphi$, we will denote by
$I_{pq}(\varphi)$ the integral
\[
I_{pq}(\varphi) =  \displaystyle \int_M \left(- p \theta_X + q
\Delta \theta_X + q {\rm Ric}(\omega_{\varphi}) + p
\omega_{\varphi} \right)^{n+1}.
\]

We need to show that
\[ {\p I_{pq} \over {\p t}}
(\varphi_t) =0.
\]

Put $\omega_t = \omega + \sqrt{-1} \p\b\p \varphi_t$. Then
\[
\begin{array}{lcl}
i_X \omega_t &=& \sqrt{-1} \b\p \left(\theta_X +
X(\varphi_t)\right)\\
i_X {\rm Ric}(\omega_t) &=& - \sqrt{-1} \b\p \Delta_t \left(
\theta_X + X(\varphi_t)\right),
\end{array}
\]
where $\Delta_t$ is the Laplacian of $\omega_t$. For simplicity,
we denote by $\Psi_t$ the function
\[
- p\;\left( \theta_X + X(\varphi_t)\right) + q\;\Delta_t \left (
\theta_X + X(\varphi_t)\right).
\]

Define
\[
\alpha_t = -p \;\p\left({\p \varphi_t \over \p t}\right) + q \;\p
\Delta_t \left({\p \varphi_t \over \p t}\right).
\]

Using the identity
\[
{\rm Ric}(\omega_t) = {\rm Ric}(\omega) - \sqrt{-1}\;\p\b\p
\log\left ( {\omega_t^n \over \omega^n}\right),
\]
we can show
\[
\sqrt{-1}\;\b\p \alpha_t = {\p  \over \p t}\left(p \;\omega_t +
q\; {\rm Ric}(\omega_t)  \right).
\]

On the other hand, we have
\[
i_X \alpha_t = {\p \Psi_t \over \p t}.
\]

This can be seen as follows: Suppose that in local coordinates,
\[
X=X^k {\p \over \p z_k}~~~{\rm and }~~~\omega_t
=\sqrt{-1}\sum_{i,j=1}^n\; g_{i \overline{j}}\; dz_i \wedge d
\overline{z_j}. \] Then
\[
\begin{array}{lcl}
{\p \Psi_t \over \p t} &=&-p\; X\left({\p \varphi_t\over \p
t}\right) + q\; {\p \over \p t}\left(g^{i\b j} {\p^2 \over \p
z_i\p\b{z}_j} \left( \theta_X + X(\varphi_t)\right)\right),\\
i_X \alpha_t &=& -p\; X\left({\p \varphi_t\over \p t}\right) +
q\; X \left( g^{i\b j}{\p^2 \over \p z_i\p\b{z}_j} \left(
{\p\varphi_t\over \p t} \right)\right).
\end{array}
\]

Notice that $\sqrt{-1} \b\p \left(\theta_X + X(\varphi_t)\right)
= i_X \omega_t$. We then have
\[
\begin{array}{lcl}
&&{\p \over \p t}\left(g^{i\b j} {\p^2 \over \p z_i\p\b{z}_j}
\left( \theta_X + X(\varphi_t)\right)\right)\\
&=& g^{i\b j} {\p^2 \over \p z_i\p\b{z}_j} \left( X^k {\p \over
\p z_k}\left( {\p \varphi_t\over \p t} \right)\right) + {\p
g^{i\b j}\over \p t} {\p
\over \p z_i} \left( X^k g_{k\b j}\right)\\
&=&g^{i\b j} {\p \over \p z_i} \left( X^k {\p^2 \over \p
z_k\p\b{z}_j}\left( {\p \varphi_t\over \p t} \right)\right) -
g^{i\b {l}}  {\p^2 \over \p z_m \p\b{z}_l}\left({\p \varphi_t
\over \p t}\right) g^{m\b{j}} {\p \over \p z_i} \left( X^k g_{k\b
j}\right)\\
&=&g^{i\b j} X^k {\p^3 \over \p z_i\p z_k\p\b{z}_j}\left( {\p
\varphi_t\over \p t} \right) - g^{i\b {l}}  {\p^2 \over \p z_m
\p\b{z}_l}\left({\p \varphi_t \over \p t}\right) g^{m\b{j}} X^k
{\p g_{k\b j} \over \p z_i}\\
&=&X^k {\p \over \p z_k}  \left( g^{i\b j}{\p^2 \over \p
z_i\p\b{z}_j} \left( {\p\varphi_t\over \p t}\right)\right).
\end{array}
\]

It follows that $i_X \alpha_t = {\p \Psi_t \over \p t}$.

For simplicity, we will denote by $R_t$ the curvature form $p\;
\omega_t + q \;{\rm Ric}(\omega_t)$. Then
\[
\sqrt{-1}\; \b \p \Psi_t ~=~ -\;i_X R_t ~~~{\rm and}~~~
\sqrt{-1}\; \b \p \alpha_t ~=~ {\p R_t \over \p t}.
\]

Hence, we have
$$
\begin{array}{lcl}
& & \sqrt{-1}\;{{\p } \over {\p t}}I_{pq}(\varphi_t) \\
& = & \int
_M \left( \sqrt{-1}\;{\p \Psi_t \over \p t} + {\p R_t \over \p t}
\right)
\left(\sqrt{-1}\;\Psi_t + R_t\right )^n \\
&=& \int _M  \left(\sqrt{-1}\;i_X \a _t + \sqrt{-1}\;\b \p \a _t
\right)\left( \sqrt{-1}\;\Psi_t + R_t\right )^n \\
&=&\int _M \sqrt{-1}\;i_X \a _t \; \left(\sqrt{-1}\;\Psi _t  +R_t\right )^n \\
& & +\; n \int _M \sqrt{-1}\;\a _t\wedge \b\p\left
(\sqrt{-1}\;\Psi_t +R_t\right )\wedge \left(\sqrt{-1}\;\Psi_t
+R_t\right)^{n-1}\\
&=&\int _M \sqrt{-1}\;i_X \a _t \; \left(\sqrt{-1}\;\Psi _t  +R_t\right )^n \\
& & - \;n \int _M \sqrt{-1}\;\a _t\wedge i_X \left
(\sqrt{-1}\;\Psi_t +R_t\right )\wedge \left(\sqrt{-1}\;\Psi_t
+R_t\right)^{n-1}  \\
& = & \int _M i_X \left(\sqrt{-1}\;\a _t\wedge\left(
\sqrt{-1}\;\Psi _t +R_t \right)^n\right).
\end{array} $$

Here we have used the second Bianchi identity: $\b \p R(g_t)=0$.
We also have used $\b\p \psi _{X,t} = -i_X R(g_t) =-i_X
(\psi_{X,t}+ R(g_t ))$.

Put
$$\eta = \sqrt{-1}\;\a _t \wedge \left(
\sqrt{-1}\;\Psi _t +R_t \right)^n
$$

We write it as $\g _0 +\cdots + \g _ {2n}$ and $i_X \eta = \beta
_0 +\beta _1 +\cdots + \beta_{2n}$. The only term which
contributes to the above integral is the $\beta _{2n}$, but
$\beta _{2n} = i_X \g _{2n+1}$ and $\g _{2n+1} =0$. Therefore,
the above integral is zero. Thus, the theorem is proved.

\section{Modified K\"ahler Ricci flow}

In the first subsection, we want to modify the K\"ahler Ricci
flow by automorphisms so that the evolved K\"ahler form is
centrally positioned with respect to a fixed K\"ahler-Einstein
metric (see Definition 6.1 below). Our argument here essentially
due to S. Bando and T. Mabuchi \cite{Bando87}.  In the second
subsection, we use this and Tian's inequality \cite{tian98} to
derive a uniform lower bound on $E_k.\; $ That in turn implies
the desired integral estimate on curvature (Corollary 6.9).
\subsection{Modified K\"ahler form by automorphisms}
  As before, let $\omega_1$  be a
K\"ahler-Einstein metric in $M$ such that  ${\rm Ric}(\omega_1) =
\omega_1.\;$ Let us first introduce the definition of "centrally
positioned":

\begin{defi}  Any K\"ahler form  $ \omega_{\varphi}$  is called centrally positioned
with respect to some K\"ahler-Einstein metric $\omega_{\rho} =
\omega + \sqrt{-1} \partial \overline{\partial} \rho$
 if it satisfies the following:
\begin{equation}
\displaystyle \int_M (\varphi -\rho) \; \theta\; {\omega_{\rho}}^n
= 0, \qquad \forall \;\theta \in \Lambda_1(\omega_{\rho}).
\label{eq:propposition}
\end{equation}
\end{defi}

We now introduce a well known functional in K\"ahler geometry
first:
\[
I( \omega_{\varphi}, \omega) = {1\over V}\displaystyle \int_M
\varphi(\omega^n - {\omega_{\varphi}}^n).
\]

 Note that this definition is symmetric with respect to $\omega$
and $\omega_{\varphi}.\;$  Alternatively, for any path
$\varphi(t) \in {\cal P}(M, \omega), $ we have
\[
  {{d\,I} \over {d\,t}} = {1\over V} \displaystyle \int_M {{\partial \varphi} \over {\partial t}}
   (\omega^n - {\omega_{\varphi}}^n)
  - {1\over V} \displaystyle \int_M \varphi \triangle_{\varphi} {{\partial \varphi} \over {\partial t}}
  {\omega_{\varphi}}^n.
\]

Put
\[
J(\omega_{\varphi},\omega) = J_{\omega}(\varphi).
\]

This implies that
\begin{equation}
{{d\,(I - J)} \over {d\,t}}  = - {1\over V} \displaystyle \int_M
\varphi \triangle_{\varphi} {{\partial \varphi} \over {\partial
t}} {\omega_{\varphi}}^n . \label{eq:Psiderivative}
\end{equation}

Now we consider a functional $\Psi$ on ${\rm Aut}_r(M)$ by,
\begin{equation}
\Psi(\sigma) = (I - J) (\omega_{\varphi}, \sigma^* \omega_1) =
(I-J)(\omega_{\varphi}, \omega_{\rho})
\label{eq:modifiedfunctional}
\end{equation}
for any $\sigma \in {\rm Aut}_r(M) $ and $\sigma^* \omega_1
=\omega_{\rho} = \omega + \sqrt{-1} \partial \overline{\partial}
\rho.\;$ If $\sigma  $ is a critical point in $\in {\rm
Aut}_r(M),\;$ then $\omega_{\rho}$ is the desired
K\"ahler-Einstein metric.
\begin{prop} Let $\omega_{\rho}$ be the minimal point of $\Psi.\;$ For any
 $u \in \Lambda_1(\omega_\rho)$, we have
 \[
   \displaystyle \int_M\; (\rho-\varphi) u  \; \omega_{\rho}^n = 0,
 \]
 or equivalently
 \[
 \rho-\varphi \perp \Lambda_1 (\omega_{\rho}).
 \]
 In other words, $\omega_{\varphi}$ is centrally positioned with
 respect to
 $\omega_{\rho}.\;$
 \end{prop}

Note that if $\Lambda_1(M) = \emptyset,$ then this proposition
hold trivially.  Before we prove this proposition, we pause to
establish the equivalence relation between  the first eigenspace
of $\omega_1$(or any K\"ahler-Einstein metric) and the space of
holomorphic vector fields (denoted by $\e (M)$).

\begin{lem}
The first eigenvalue of $\Delta _{\omega_1} \geq 1.\;$ Moreover,
there is a 1-1 correspondence between the first eigenspace
$\Lambda_1$ of $\omega_1 $ and the space of holomorphic vector
fields $\e (M).\;$
\end{lem}
The lemma is well-known. For the reader's convenience, we outline
its proof here.
\begin{proof}  Let $\lambda_1$ be the first eigenvalues of
$\omega_1$  and  $u$ is any eigenfunction  of $ \omega_1$ with
eigenvalue $ \lambda_1, \;$ so  $\Delta _{\omega_1} u = -\lambda_1
u.\;$  Define a vector field $X$ by $i_X \omega_1 = \sqrt{-1}
\bar \p \;u.\;$ Then by a direct computation, we have
 \[
\displaystyle \int_M\; \mid \bar \p \;X\mid^2\; {\omega_1}^n =
\lambda_1^2 \displaystyle \int_M\; u^2\; {\omega_1}^n -
\displaystyle \int_M\; \mid \p \;u \mid^2\; {\omega_1}^n.
  \]
  This implies that

$$
\begin{array}{lcl}
 \lambda_1^2\;\int _M u^2 \; {\omega_1}^n &  = & \int _M | \p u | ^2\; {\omega_1}^n
        + \displaystyle \int_M\; \mid \bar \p \;X\mid^2\;
        {\omega_1}^n\\
       & \ge &  \lambda_1 \int _M u^2\; {\omega_1}^n + \displaystyle \int_M\; \mid \bar \p \;X\mid^2\;
        {\omega_1}^n \\
        & \ge & \lambda_1 \int _M u^2\; {\omega_1}^n.
        \end{array}
$$
Here we have used the variational characterization of $\lambda
_1.\;$ Thus $\lambda_1 \geq 1.\; $ If the equality holds, i.e.,
$\lambda_1 = 1$, we  have
 $\bar \p X  = 0.\;$ It follows that $X$ is a holomorphic vector
 field.

Conversely, if $X$ is a holomorphic vector field, we define $u$ by
 $i _X \omega_1 = \b \p u $ and
 $ \displaystyle \int_M \; u\; {\omega_1}^n =0.\;$ Then a straightforward computation shows that
 \[
 \b \p ( \triangle_{\omega_1}\; u + u) = 0.
 \]
 It follows that $u$ is an eigenfunction with
eigenvalue $1$.  So we have established the following
identification
$$
\eta (M) \simeq \{ {\rm \, \, eigenfunctions \, \, of \, \, } \o
_1 {\rm \, \, with \, \, eigenvalue \, \,} 1 \}.
$$
\end{proof}
Next we return to the proof of proposition 6.2.

\begin{proof} Let $\sigma_s$ be the one parameter subgroup
generated by the real part of
 $\overline{\partial} u$, write
 \[
 \begin{array}{lcl}
 \omega_{\rho_s} & = & \sigma_s^* \omega_{\rho} = \omega_{\rho} + \sqrt{-1} \partial \overline{\partial}
 \left(\rho_s  -\rho\right)  \\& = & \omega_{\varphi} + \sqrt{-1} \partial \overline{\partial} ( \rho_s -\varphi).
 \end{array}\]

 One can easily see that ${{d\, \rho_s} \over {d\,s}}\mid_{s=0} = u $ modulo constants.  Denote
 the complex Laplacian operator of $\omega_{\rho}$ by $\triangle_{\rho}.\;$ Then,
 \[
     \triangle_{\rho} u + u = 0, \qquad \forall \; u \in \Lambda_1(\omega_{\rho}).
 \]

 Computing the derivative of  $\Psi$ (see Formula \ref{eq:Psiderivative}) along this holomorphic path,
  we have
 \[
\begin{array}{lcl} 0 &=&  {d \over {d\,s}} \Psi(\sigma_s)\mid_{s=0} \\
& = & - { 1 \over V}\;\displaystyle \int_M\; (\rho-\varphi)\;
\triangle_{\rho}\;{{d\, \rho_s} \over {d\,s}}\mid_{s=0} \;
\omega_{\rho}^n
\\ & = & - { 1 \over V}\;\displaystyle \int_M\; (\rho-\varphi)\;\triangle_{\rho}\; u\; \omega_{\rho}^n
\\ & = & { 1 \over V}\;\displaystyle \int_M\; (\rho-\varphi)\; u\; \omega_{\rho}^n.
\end{array}
 \]

 In other words,
 \[
   \rho-\varphi \perp \Lambda_1 (\omega_{\rho}).
 \]
\end{proof}
The rest of the subsection is devoted to prove that there always
exists a minimizer of $\Psi$ in ${\rm Aut}_r(M).\;$
 Recall that $\omega_1$ is a
K\"ahler-Einstein metric, so $\omega_{\rho}$ is also
K\"ahler-Einstein metric:
\begin{equation}
 \left(\omega_{\varphi} + \sqrt{-1} \p \bar \p (\rho(t)-\varphi(t)) \right)^n = {\omega_{\rho}}^n = e^{-(\rho - \varphi) + h_{\varphi} } {\omega_{\varphi}}^n,
\label{eq:ke}
\end{equation}
where
\[
{\rm Ric}(\omega_{\varphi}) - \omega_{\varphi} = \sqrt{-1}
 \partial \overline{\partial} h_{\varphi}.
\]

We shall normalize $h_{\varphi}$ and $\rho$ as
\[
  {1\over V} \displaystyle  \int_M \; e^{-(\rho - \varphi) + h_{\varphi} } \omega_{\varphi}^n =
   {1\over V} \displaystyle  \int_M \; e^{h_{\varphi} } \omega_{\varphi}^n = 1 .
\]

Therefore, we have
\[
  \displaystyle \sup_M (\rho-\varphi) \geq 0.
\]

\begin{prop} The following inequalities hold
\begin{equation}
 I - J \leq I \leq (n+1)(I - J).
 \label{eq:IJcompare}
\end{equation}
\end{prop}

\begin{proof} From the definition
of $I,$ we have
\[
\begin{array}{lcl} I(\omega_{\varphi}, \omega) & = &
 {1\over V} \displaystyle \int_M \varphi(\omega^n - {\omega_{\varphi}}^n) \\
 & = & {1\over V}  \displaystyle \int_M \varphi \wedge (- \sqrt{-1} \partial \overline{\partial} \varphi)
  \displaystyle \sum_{i=0}^{n-1} \omega^i \wedge \omega_{\varphi}^{n-i-1}\\
  & = & \displaystyle \sum_{i=0}^{n-1} {1\over V} \displaystyle \int_M  \;\sqrt{-1}\;\partial \varphi
  \wedge \overline{\partial} \varphi
  \wedge \omega^i \wedge \omega_{\varphi}^{n-i-1} \geq 0.
\end{array}
\]

Therefore,
\[
  I -J = \displaystyle \sum_{i=0}^{n-1} {{n-i} \over {n+1}} {1\over V} \displaystyle \int_M \;\sqrt{-1}\;\partial \varphi
  \wedge \overline{\partial} \varphi
  \wedge \omega^i \wedge \omega_{\varphi}^{n-i-1}.
\]

This in turn implies that
\[
 I - J \leq I \leq (n+1)(I - J).
\]
\end{proof}

Next we can prove that $\Psi$ always achieve its minimum value in
${\rm Aut}_r(M).\;$

\begin{lem} The minimal value of $\Psi$ can be attained in ${\rm Aut}_r(M).\;$ Moreover, $\Psi$
is proper.
\end{lem}

\begin{proof} Observe that
\[
\Psi(\sigma) = (I - J)(\omega_{\varphi},\omega_{\rho}) \geq {
1\over {n+1}} \int_M (\rho -\varphi)(\omega_{\varphi}^n
              - \omega_{\rho}^n)\geq 0.
\]

Put $G_r = \{\sigma \in Aut_r(M)| \Psi(\sigma) \leq r\}$ and $E_r
= \{\rho | \sigma^* \omega_1
    = \omega_{\rho}, \sigma \in G_r\}.\;$
Then,
\[
\int_M \; (\rho -\varphi) (\omega_{\varphi}^n - \omega_{\rho}^n
)\leq {(n+1) r }, \qquad \forall \rho \in E_r.
\]

However,
\[
\begin{array}{lll}
 - {1\over V} \displaystyle \int_M \; (\rho - \varphi) \omega_{\rho}^n & = & -{1\over V} \displaystyle \int_M \; (\rho - \varphi) e^{-(\rho - \varphi)
  + h_{\varphi} } \omega_{\varphi}^n
\\
 & \geq & - C_2 \int_M \; \rho\; e^{-\rho}\;\omega_{\varphi}^n
 -C_2'\\
 &\geq & - C_3.
 \end{array}  \]

 Therefore, we have
\[
   \int_M\; (\rho-\varphi) \omega_{\varphi}^n \leq C_3'.
\]

  Since $\triangle_{\varphi} (\rho-\varphi) \geq -n,\;$ by the Green formula, we have
  \[
  \begin{array}{lcl} & &
   \displaystyle \sup_M (\rho-\varphi) \\ & \leq & {1\over V} \displaystyle \int_M (\rho- \varphi) \omega_{\varphi}^n
   - \displaystyle \max_{x \in M}\;  \left( {1\over V} \int_M (G(x,\cdot) + C_4)
   \triangle_{\varphi}(\rho-\varphi)
   \omega_{\varphi}^n (y)\right)\\
   & \leq & {1\over V} \displaystyle \int_M (\rho-\varphi) \omega_{\varphi}^n + n C_4,
  \end{array}
  \]
  where $G(x,y)$ is the Green function associated to $\omega_{\varphi}$ satisfying
  $G(x,\cdot)\geq 0.\;$ Therefore, there exists a uniform constant $C$ such that
  \[
  \displaystyle \sup_{M}\; (\rho-\varphi) \leq C.
  \]

  On the other hand, we have
  \[
  \begin{array}{lcl} -\int_M (\rho-\varphi) \omega_{\rho}^n & \leq & (n+1) r
        - \int_M (\rho-\varphi) \omega_{\varphi}^n
  \\ & \leq & (n+1) r  + n C_5 - \displaystyle \sup_M (\rho-\varphi)
  \\ & \leq & C_6.
  \end{array}
  \]

  Following Proposition 6.6 below, we can prove that there exists a constant $C $ such that

    \[
                \displaystyle \inf_M \; (\rho-\varphi) \geq - C.
    \]

    Hence, by the $C^2$ estimate of Yau \cite{Yau78} and $C^3$ estimate of Calabi, we obtain
    \[
     \|\rho-\varphi\|_{C^3} \leq C_7(r), \qquad \;\forall \rho \in E_r.
    \]

  Then $E_r$ is compact in $C^2$ topology, and so is $G_r.\;$ In particular,
  the minimal value of $\Psi$ can be attained.
  \end{proof}

  \begin{prop} Let $\omega_{\rho}$ be a K\"ahler-Einstein metric,
  then
$$
0 \le -\inf _M (\rho-\varphi) \le C \Bl {1 \over V} \int _M
(-(\rho-\varphi) ) \o _{\rho} ^n +1 \Br.
$$
\end{prop}

The proposition is known (cf. \cite{tian98}), we include its
proof here for reader's convenience.

\begin{proof} Denote by $\Delta _{\rho}$ the Laplacian of $\o
_{\rho}.\;$ Then, because $\o_{\varphi} + \p \b \p (\rho-\varphi)
>0$, we see that $\o_{\varphi} = \o _{\rho} - \p \b \p
(\rho-\varphi)
>0$ and taking the trace of this latter expression with respect
to $\o _{\rho}$, we get
$$
n - \Delta _{\rho} (\rho-\varphi) = \tr _{\o _{\rho}} \o_{\varphi}
> 0.
$$

Defining now $(\rho-\varphi) _- (x) = \max \{ -(\rho-\varphi) (x)
, 1 \} \ge 1$, so that
 \[
 {(\rho-\varphi) _-} ^p (n - \Delta
_{\rho} (\rho-\varphi) ) \ge 0. \]
 And integrating this, we
get
\begin{align*}
0 & \le {1 \over V} \int _M {(\rho-\varphi) _- }^p (n - \Delta _{\rho} (\rho-\varphi) ) \o _{\rho} ^n \\
  &= {n \over V} \int _M {(\rho-\varphi) _- }^p \o _{\rho} ^n + {1 \over V} \int _M
     \nabla _{\rho} {(\rho-\varphi) _- }^p \nabla _{\rho} (\rho-\varphi) \o _{\rho} ^n \\
  &= {n \over V} \int _M {(\rho-\varphi) _- }^p \o _{\rho} ^n + {1 \over V} \displaystyle \int _{\{ (\rho-\varphi)
     \le -1 \} } \nabla _{\rho} {(\rho-\varphi) _- }^p \nabla _{\rho} (\rho-\varphi) \o _{\rho} ^n \\
  &= {n \over V} \int _M {(\rho-\varphi) _- }^p \o _{\rho} ^n + {1 \over V} \int _M
      \nabla _{\rho} {(\rho-\varphi) _- }^p \nabla _{\rho} (-(\rho-\varphi) _- ) \o _{\rho} ^n \\
  &= {n \over V} \int _M {(\rho-\varphi) _- }^p \o _{\rho} ^n  - {1 \over V} {4p \over
     {(p+1 )^2}} \int _M | \nabla _{\rho} {(\rho-\varphi) _- }^{{{p+1} \over 2}} |^2
      \o _{\rho} ^n,
\end{align*}
which yields, using the fact that $(\rho-\varphi) _- \ge 1$ and
hence ${(\rho-\varphi) _- }^p \le {(\rho-\varphi) _- }^{p+1} ,$
$$
{1 \over V} \int _M \bigl| \nabla _{\rho} {(\rho-\varphi) _-}
^{{{p+1} \over 2}} \bigr|^2 \o _{\rho} ^n  \le {{n(p+1)^2 } \over
{4p V}} \int _M {(\rho-\varphi) _- }^{p+1} \o _{\rho} ^n .
$$


Note that $\omega_{\rho}$ is a K\"ahler Einstein metric which has
a uniform Sobolev constant. Thus, we have
$$
{1 \over V} \Bl\int_M  |(\rho-\varphi) _- | ^{{(p+1) n} \over
{n-1}} \o _{\rho} ^n \Br ^{{n-1} \over n} \le {{c(p+1)} \over V}
\int _M {(\rho-\varphi) _- }^{p+1} \o _{\rho} ^n.
$$

Moser's iteration will show us that
$$
\sup _M (\rho-\varphi) _- = \lim _{p \to \infty} \|
(\rho-\varphi) _- \| _{ L^{p+1} (M,\o _{\rho})} \le C \|
(\rho-\varphi) _- \| _{L^2 (M, \o _{\rho} )}.
$$

Recall that $\l _1 (\o _{\rho} ) \ge 1$, so that the Poincar\'e
inequality reads
\begin{align*}
{1 \over V} \int _M \Bl (\rho-\varphi) _- - {1 \over V} \int _M
(\rho-\varphi) _- \o _{\rho} ^n \Br ^2 \o _{\rho} ^n
&\le {1 \over V} \int _M | \nabla (\rho-\varphi) _- | ^2 \o _{\rho} ^n \\
&\le {C \over V} \int _M (\rho-\varphi) _- \o _{\rho} ^n,
\end{align*}
where we have set $p=1$ and used the same reasoning as before.
This then implies that
$$
\max  \{ - \inf _M (\rho-\varphi) , 1 \} = \sup _M (\rho-\varphi)
_- \le {{C} \over V} \int _M (\rho-\varphi) _- \o _{\rho} ^n,
$$
since $\int_M e^{-h_\varphi + (\rho-\varphi)} \o_{\rho}^n = V$,
we can easily deduce $\displaystyle \int_{(\rho-\varphi) > 0}
(\rho-\varphi) \o_{\rho}^n \le C$. Combining this together with
the above, we get
$$
-\inf _M (\rho-\varphi) \le {C \over V} \int _M (-(\rho-\varphi)
) \o _{\rho} ^n + C,
$$
which proves the proposition.

\end{proof}

\subsection{Application to the K\"ahler Ricci flow}

Let $\varphi(t)$ be the global solution of the K\"ahler Ricci
flow in the level of K\"ahler potentials. According to Lemma 6.5,
there exists a one parameter family of K\"ahler Einstein metrics
$\omega_{\rho(t)} = \omega + \sqrt{-1} \partial \bar \partial
\rho(t) $ such that $\omega_{\varphi(t)}$ is centrally positioned
with respect to  $\omega_{\rho(t)} $ for any $t \geq 0. \;$
Suppose that $\omega_{\varphi(0)}$ is already centrally
positioned with the K\"ahler-Einstein metric $\omega_1 = \omega +
\sqrt{-1} \p \bar \p \rho(0).\;$ Recall that
$E_{k,\omega}(\varphi)$ and $\nu_{\omega}$ all satisfy the cocycle
condition:
\[
E_{k,\omega} (\varphi) +  E_{k,\omega_\varphi}(\psi -\varphi) =
E_{k,\omega}(\psi)
\]
for any $k=0,1,\cdots, n.\;$ Note that $\nu_{\omega} =
E_{0,\omega}.\;$

\begin{theo} On a K\"ahler-Einstein manifold, the $K$-energy $\nu_{\omega}$ is uniformly
bounded from above and below along the K\"ahler Ricci flow.
Moreover, there exist some uniform constants $c,C,\;C'$ and $C''$
such that
\[
\begin{array}{ccl}
|J_{k,\omega_{\rho(t)}}(\varphi(t) - \rho(t) )| & \leq &
\{\nu_{\omega}(\varphi(t)) + C\}^{1\over \delta },\\
\log { {\omega_{\varphi}}^n \over {{\omega_{\rho(t)}}^n}} & \geq &
 - 4 C''\, e^{2\left(\nu_{\omega}(\varphi(t)) + C\right)^{1\over \delta} + C')},\\
E_k (\varphi(t)) & \geq &  - e^{c\left(1 + {\rm max} \{0,
\nu_\omega(\varphi(t))\} + \left(\nu_{\omega}(\varphi(t)) +
C\right)^{1\over \delta }\right) }.
\end{array}
\]
\end{theo}
\begin{proof}
 Since $\omega_{\varphi(t)}$ is centrally positioned
with respect to the K\"ahler-Einstein metric $\omega_{\rho(t)},$
Proposition 6.1 implies that
\[ \varphi(t) - \rho(t) \perp \Lambda_1 (\omega_{\rho(t)}).
 \]

  Theorem 3.1  implies
that K-energy is proper with respect to the evolved K\"ahler
metric $\omega_{\varphi}$ and the modified K\"ahler-Einstein
metric $\omega_{\rho(t)}.\;$ Thus,
\[
  \nu_{\omega_{\rho(t)}}(\varphi(t) - \rho(t) ) \geq \left(J_{\omega_{\rho(t)}}(\varphi - \rho(t))\right)^{\delta} - c
\]
for some uniform constant $\delta > 0$ and $c.\;$  Since the K
energy satisfies the cocycle condition, we have
\[ \nu_{\omega}(\varphi(t)) - \nu_{\omega_{\rho(t)}} (\varphi(t)
-\rho(t)) = \nu_{\omega} (\rho(t)).
\]

Lemma 3.8 implies that the K-energy monotonely decreases along the
K\"ahler Ricci flow
 \[
\nu_{\omega}(\varphi(t)) \leq \nu_{\omega}(\varphi(0)), \qquad
\forall\; t \; < \infty.
  \]

Combining the three inequalities above, we arrive at
\[
\nu_{\omega}(\varphi(0)) \geq \nu_{\omega}(\varphi(t)) \geq
\left(J_{\omega_{\rho(t)}}(\varphi - \rho(t))\right)^{\delta} - c
+ \nu_{\omega} (\rho(t)).
\]

Note that the K energy is invariant under automorphisms and the
fact that $\omega_{\rho(t)}$ is path connected with $\omega_1$
via automorphisms, then we have
\[
\nu_{\omega} (\rho(t)) = \nu_{\omega}(\rho(0)).
\]
Thus
\[
0\leq J_{\omega_{\rho(t)}}(\varphi - \rho(t)) \leq
\left(\nu_{\omega}(\varphi(t)) + C \right)^{1 \over \delta} \leq
\left(\nu_{\omega}(\varphi(0)) + C \right)^{1 \over \delta}.
\]

In particular, the K energy has a uniform up-bound and lower
bound along the K\"ahler Ricci flow. Lemma 4.10 implies that
$E_{k,\omega_{\rho}} (\varphi(t) - \rho(t)) $ are uniformly
bounded from below. Now,

\[
E_{k,\omega} (\varphi(t)) =  E_{k,\omega_\rho}(\varphi(t) -
\rho(t) ) + E_{k,\omega}(\rho(t))
\]

Similarly since $E_k$ is invariant under automorphisms, we have
\[
E_{k,\omega}(\rho(t)) = E_{k,\omega}(\rho(0)).
\]

Thus
\[
\begin{array}{lcl}
E_{k,\omega} (\varphi(t)) & = & E_{k,\omega_\rho}(\rho(t)-\varphi(t)) + E_{k,\omega}(\rho(0)) \\
& \geq &   - e^{c\left(1 + {\rm max} \{0, \nu_\omega(\varphi(t))\}
+ J_{\omega_\rho}(\varphi - \rho)\right)} + E_{k,\omega}(\rho(0))
\\& = & - e^{c\left(1 + {\rm max} \{0, \nu_\omega(\varphi(t))\} +
\left(\nu_{\omega}(\varphi(t)) + C\right)^{1 \over \delta}\right)
} - C_1,
\end{array}
\]
where $c, C$ and $ C_1$ are some uniform constant.  It also
implies that (from the explicit expression of the K energy
(\ref{eq:kenergyexp0})):
\[
  \displaystyle \int_M \; \left(\ln {{\omega_{{\varphi}(t)}}^n \over
  {\omega_{\rho}}^n } \right) {\omega_{{\varphi_{\rho}}(t)}}^n \leq
  \left(\nu_{\omega}(\varphi(t)) + C\right)^{1\over \delta} + C_2,
\]
where $C_2$ is some uniform constant.  Proposition 3.13 then
implies that $\log { {\omega_{\varphi}}^n \over
{{\omega_{\rho(t)}}^n}} $ is uniformly bounded from below:
\[
\begin{array}{lcl}
\displaystyle \inf_{M}\; \left( \log {\omega_{\varphi}^n \over
{\omega_{\rho}}^n } \right)(x)& \geq &- 4 C_3\, e^{2(1+
\displaystyle \int_M \; \left( \log {\omega_{\varphi}^n \over
{\omega_{\rho}}^n }
\right)\;{\omega_{\rho}}^n )} \\
& \geq & - 4 C_3\, e^{2\left(\nu_{\omega}(\varphi(t)) +
C\right)^{1 \over \delta} + C')}

\end{array}
\]
where $C,C_3 $ and $ C'$ are some uniform constant.  Corollary 4.5
shows that $J_k (k=0,1,\cdots n-2) $ are uniformly bounded from
above and below.
\end{proof}
An immediate corollary is
\begin{cor} The energy functional $E_k (k=0,1,\cdots n) $ has a
uniform lower bound from below along the K\"ahler Ricci flow.
 \end{cor}
\begin{proof} Since $E_k$ is invariant under action of
automorphisms. Thus
\[
E_k(\varphi) = E_k(\tilde{\varphi}) \geq -C.
\]
\end{proof}
Now combing Theorem 4.10 and this Corollary, we arrive at the
following important corollary:

\begin{cor} For each $k=0,1,\cdots n,$ there exists a uniform constant $C$
such that the following holds (for any $T\leq \infty$) along the
K\"ahler Ricci flow:
\[
\displaystyle \int_{0}^{T}\; {{k+1}\over V}\; \displaystyle
\int_M\; \left(R(\omega_{\varphi(t)}) - r\right) \;  {{\rm
Ric}(\omega_{\varphi(t)})}^{k} \wedge {\omega_{\varphi(t)}}^{n-k}
\;d\,t \leq C.
\]

When $k=1, $ we have
\[
  \displaystyle \int_0^{\infty}\;{{1}\over V}\;
 \displaystyle \int_M\; (R(\omega_{\varphi(t)})-r) ^2\;   {\omega_{\varphi(t)}}^{n}  \;d\,t
 \leq C < \infty.
\]
\end{cor}

\section{Injectivity Radius}
In 1959, Klingenberg proved that for any compact oriented, even
dimensional manifold without boundary, if the sectional curvature
is bounded in $(0,1],$ then the injectivity radius  is at least
$\pi.\;$  This theorem of Klingenberg does not apply to the
evolved metrics in the K\"ahler Ricci flow since we do not know
if the positivity of the sectional curvature will be  preserved.
However, by  Theorem 2.2, the bisectional curvature is positive
along the K\"ahler Ricci flow if the initial metric has a
positive bisectional curvature. Therefore, we need to adopt
Klingenberg's original theorem to our case. Namely, obtaining a
similar estimate of the injectivity radius based on the
positivity of the bisectional curvature only. Such a lemma is a
natural extension of the original Klingenberg's theorem to the
K\"ahler setting.

\begin{lem} Suppose that $(M,g)$ is an orientable compact
K\"ahler manifold with bisectional curvature bounded in
$(0,1].\;$ Then there exists some uniform constant $\beta > 0$
such that the injectivity radius must be no less than $\beta
\pi\;$\footnote{According to Corollary 12.2, the best constant is
$\beta =  { 1 \over \sqrt{2}}.\;\;$}.
\end{lem}
\begin{proof} We  follow the arguments in the proof of  Klingenberg
theorem (c.f.\cite{CE75}). Since the bisectional curvature $\leq
1,$ there exists a uniform constant ${1 \over {\beta}^2} $ such
that the sectional
 curvature is  uniformly bounded from above by ${1\over {\beta}^2}.\;$ This follows that
  the conjugate radius is not shorter than
$\beta \pi:$
\[
conj\; rad_M \geq \beta \pi.
\]

A lemma in \cite{CE75} by Cheeger and Ebin  asserts:
\[
  inj_M = \displaystyle \min \{ \beta \pi, {1\over 2} {\rm the \;length\; of \;shortest \; closed \; geodesic}\}.
\]

We want to prove the lemma by contradiction.   If the injective
radius $< \beta \pi,$ then there exists a shortest closed
geodesic which realizes this injectivity radius.
 Denote this shortest closed curve by $ c_0(t)\;(0\leq t \leq 2\;inj_M)$
parameterized by the  arc length. Suppose $J$ is the underlying
complex structure. The plane spanned by $c_0(t)' $ and $J
(c_0(t)')$ is a holomorphic plane. Thus the sectional curvature of
this plane must be strictly positive. Deform $c_0$ on the
direction of $J(c_0(t)').\;$ Since $J(c_0(t)')$ is a parallel
vector field along this closed geodesic, the second variation
 in this direction is strictly negative.
 Therefore, there exists a 1-parameter family of nearby closed
 curves
$c_s:\mathbb{R} /Z \rightarrow V, t \rightarrow exp_{c_0(t)} (s
J(c_0(t)'))\;$ which are strictly shorter than $c_0$ provided $s$
is small enough. Since the length of $c_0$ equals to $2 \;
inj_M,$   the entire curve $c_s \;(s>0)$ must be contained in the
closed ball with radius $ \leq {1\over 2} L(c_s) < inj_M.\;$
Thus, one can lift the entire curve $c_s(t)$ as a closed curve
$\tilde{c}_s$ in $T_{c_s(0)}M$ such that $\tilde{c}_s(0) = 0.\;$
Since everything occurs within the conjugate radius, by taking
limit, we can  lift up $c_0(t)$ as a closed curve in
$T_{c_0(0)}M.\;$ That is a contradiction since the lifting of
$c_0$ is a straight line.
\end{proof}

\section{Harnack inequality}
Recall Cao's Harnack inequality in the K\"ahler Ricci flow:
\begin{theo} \cite{Cao92} Let $g_{i\overline{j}}$ be the solution of the K\"ahler Ricci flow with positive bisectional curvature. Then
for any $x,y \in M$ and $0< t_1 <t_2 <\infty,$  the scalar
curvature $R$ satisfies the inequality:
\[
    R(x,t_1) \leq {{e^{t_2} -1}\over {e^{t_1}-1}} e^{\Delta \over 4} R(y,t_2).
\]
Here $\Delta$ is defined as
\[
  \Delta = \Delta(x,y,t_1,t_2) = \displaystyle \inf_{\gamma} \displaystyle \int_{t_1}^{t_2}\;
|{\gamma}^{'}(s)| ^2\;d\,s
\]
where the infimum is taken over all curves from $x$ to $y$, where $|{\gamma}^{'}(s)|_s$  is the
velocity of $\gamma$ at time $s$.
\end{theo}
The basic ideas\footnote{These types of arguments are due to R.
Hamilton in the real case (cf. \cite{Hamilton93}). } of the proof
in \cite{Cao92} can be described as follows: If $g$ is a K\"ahler
Ricci soliton, we have
\begin{equation}
  R_{i \overline{j}} - g_{i \overline{j}} = f_{,i \overline{j}}
\label{eq:solition}
\end{equation}
and
\[
  f_{,i j}=0,\qquad \forall\; i,j=1,2,\cdots, n.
\]

 Thus, $ X = f^{,i}{{\partial }\over {\partial z^i}} $ is a
holomorphic vector field. Taking Laplacian of the soliton
equation (\ref{eq:solition}),
 we arrive at the following
\[
\Delta R_{i \overline{j}} + R_{i \overline{j} k \overline{l}} R_{l \overline{k}} - R_{i \overline{j}} + R_{i \overline{j},k} \,f_{,\overline{k}} + R_{i \overline{j},\overline{k}} f_{,k} + R_{i \overline{j} k \overline{l}}\, f_{,l} f_{,\overline{k}} = 0.
\]

Motivated  by this identity for Ricci solitons,
  Cao  introduced the following 2-tensor $Q_{i \overline{j}}$
 for any vector $v\in T_x V,$
\begin{eqnarray}
Q_{i\overline{j}} & = & \Delta R_{i \overline{j}} + R_{i
\overline{j} k \overline{l}} R_{l \overline{k}} - R_{i
\overline{j}} + R_{i \overline{j},k} \,v^{,{k}} + R_{i
\overline{j},\overline{k}} v^{,\overline{k}} + R_{i \overline{j}
k \overline{l}}\, v^{,\overline{l}} v^{,{k}} + {{R_{i
\overline{j}}} \over {1 - e^{-t}}}\nonumber \\ & = &
{{\partial}\over {\partial t}} R_{i \overline{j}} + R_{i
\overline{k}} R_{k\overline{j}} \nonumber \\
& & \qquad \qquad- R_{i \overline{j}} + R_{i \overline{j},k}
\,v^{,\overline{k}} + R_{i \overline{j},\overline{k}} v^{,k} +
R_{i \overline{j} k \overline{l}}\, v^{,l} v^{,\overline{k}} +
{{R_{i \overline{j}}} \over {1 - e^{-t}}}. \label{eq:solition2}
\end{eqnarray}

Clearly, $Q$ is a positive tensor at $t=0$ and $t=\infty.\;$
Through  tedious  but direct calculations, Cao proved that $Q$ is
positive for all the time and for all vectors $v(x,t).\;$  Taking
trace on both side of (\ref{eq:solition2}), we obtain the
following
\[
{{\partial R}\over {\partial t}} + R_{,k} v^{,k} + R_{,\overline{k}} v^{,\overline{k}} + R_{i \overline{j}} v^{\overline{i}} v^{j} +  {{R} \over {1 - e^{-t}}} > 0.
\]

Let $v_k = {{-R_{,k} \over R}} $, then
\[
{{\partial R}\over {\partial t}} -{{|D\,R|^2}\over {R}} +  {{R} \over {1 - e^{-t}}}> 0.
\]

Using this inequality and a similar argument of Li-Yau
\cite{PYau86},  Cao \cite{Cao92} proved  the Harnack inequality
for the scalar curvature of $M.\;$

\section{Convergence by sequence in any $C^l$ norm}
In this section, we want to show that for any sequence of metrics
over the K\"ahler Ricci flow, there exists a subsequence which
converges to a K\"ahler-Einstein metric with constant bisectional
curvature. We first prove that the bisectional curvature and its
derivatives are uniformly bounded in complex dimension 2 in the
first subsection. In the second subsection, we then prove the
convergence by sequences.
\subsection{Uniform curvature bound in complex dimension 2} In
this subsection,we concentrate on complex dimension 2 and we will
prove that the scalar curvature is uniformly bounded from above
along the flow. One should note that only Lemma 9.3
 need to be
proved in complex dimension 2.  All other theorems, lemmas hold
for all dimensions.

\begin{lem} In the K\"ahler Ricci flow with positive bisectional curvature, denote the maximal scalar
 curvature at time $t$ as $R_{max}(t).\; $  Then
 \[
    R_{max}(t) \leq 2 R_{max(t_0)}, \qquad \forall \; t \in [t_0,
    t_0 + { 1 \over {2 R_{max}(t_0)}}].
 \]
\end{lem}
\begin{proof}  During time $t \in [t_0, t_0 + {1 \over {2
R_{max}(t_0)}}],$ we have
\[
  {d\over {d\,t}} R_{max} \leq {R_{max}}^2.
\]

Thus,
\[
  R_{max}(t) \leq 2 R_{max}(t_0), \; \forall\; t \in [t_0,t_0 + {1 \over {2
  R_{max}(t)}}].
\]
\end{proof}

 By Theorem 4.10 and Corollary 6.9, for any fixed period $T$\footnote{The value of $T$ will be
fixed later in the subsection when we prove Theorem 9.4.}, we have

\[
\begin{array}{lcl} \displaystyle \int_{0}^{\infty}  \displaystyle \int_M\; (R-r) ^2\;   \omega_{\varphi}^{n}  \;d\,t
& = & \displaystyle \sum_{n=0}^{\infty} \displaystyle \int_{nT}^{(n+1)T}\; {{1}\over V}
\; \displaystyle \int_M\; (R-r) ^2\;   \omega_{\varphi}^{n}  \;d\,t \\
& = &
\displaystyle \int_{0}^{\infty}\; {{1}\over V}\; \displaystyle \int_M\; (R-r) ^2\;   \omega_{\varphi}^{n}  \;d\,t\\
 &  = &
    \displaystyle \int_{0}^{\infty}\; {{1}\over V}\; \displaystyle \int_M\; (R-r)\;  {\rm Ric} \wedge
    \omega_{\varphi}^{n-1} \;d\,t< \infty.
\end{array}
\]

Thus,
\[
 \displaystyle \lim_{n \rightarrow \infty} \displaystyle \int_{nT}^{(n+1)T}\; {{1}\over V}\;
 \displaystyle \int_M\; (R-r) ^2\;   {\omega_{\varphi}}^{n}  \;d\,t = 0.
\]

This follows that $ \int_M \; (R-r) ^2\;
{\omega_{\varphi(t)}}^{n}$ is  small for almost all $t$ large. In
other words, at every interval of length $T,$ there exists at
least one time $t$ such that this integral is small:
\begin{lem} For any sequence $s_i \rightarrow \infty,$ and for any fixed time
period $T$, there exists $t_i \rightarrow \infty$ and $0 < s_i -t_i < T$ such that
\begin{equation}
\displaystyle \lim_{t_i \rightarrow \infty} {{1}\over V}\;
\displaystyle \int_M\; (R-r) ^2\;   \omega_{\varphi}^{n}  =0.
\label{eq:sequence}
\end{equation}
\end{lem}

\begin{lem} Over the K\"ahler Ricci flow on a K\"ahler surface, if $t_i\rightarrow \infty$ satisfies the condition (\ref{eq:sequence}),
then $R_{max}(t_i) = \displaystyle \max_{p \in V} R(p)(t_i)$ is
uniformly bounded from above.
\end{lem}

\begin{proof}
 Choose time $\tau_i < t_i$ such that $t_i -\tau_i = {1\over {2 R_{max}(\tau_i)}}$\footnote{In Hamilton's paper,
  Hamilton choose $\tau_i$ in a different way.}. Such a $\tau_i$ can always be chosen. Following Lemma 9.1, we
have
\[
R_{max}(t) \leq 2 R_{max}(\tau_i),\;\forall t\in [\tau_i,t_i].
\]

Recall the flow equation,
\[
   {{\partial}\over {\partial t}} g_{\alpha \overline{\beta}} = g_{\alpha \overline{\beta}} - R_{\alpha
   \overline{\beta}}.
\]

Thus, the distance grows at most by a constant factor since $
R_{\alpha \overline{\beta}} > 0$ for all the time.
 For any fixed point $p$, we have
\[
    \left(g_{i\overline{j}} \right)_{n\times n} (p,t) \leq  \left(g_{i\overline{j}} \right)_{n\times n} (p,\tau_i),
    \qquad\forall \;t\in [\tau_i,t_i].
\]

On the other hand,
\[
\begin{array}{lcl}  {{\partial}\over {\partial t}} g_{\alpha \overline{\beta}} & = & g_{\alpha \overline{\beta}}
- R_{\alpha \overline{\beta}}\\
& \geq & - R_{max}(t)  g_{\alpha \overline{\beta}} \\
& \geq & - 2 R_{max}(\tau_i)\; g_{\alpha \overline{\beta}},\qquad \forall\;t \in [\tau_i,t_i].
\end{array}
\]

Thus,
\[
\begin{array}{lcl}
\left(g_{i\overline{j}} \right)_{n\times n} (p,t_i) & \geq &\left(g_{i\overline{j}} \right)_{n\times n} (p,\tau_i)
\cdot e^{-2  R_{max}(\tau_i) (t_i-\tau_i)}
\\
& =&\left(g_{i\overline{j}} \right)_{n\times n} (p,\tau_i) e^{-1}.
\end{array}
\]

Therefore,
\[
  \left(g_{i\overline{j}} \right)_{n\times n} (p,t) \leq 3 \cdot \left(g_{i\overline{j}} \right)_{n\times n} (p,t_i),
  \qquad \forall\;t\in [\tau_i,t_i].
\]

If $d(\xi,X)$ is the geodesic distance at time $t_i , $ then
\[
\Delta (\xi,\tau_i,X,t_i) \leq 9 {{{d(\xi,X)}^2}\over
{t_i-\tau_i}}.
\]

For all $X$ in a ball around $\xi$ of radius
\[
             \rho = {{\pi}\over{\sqrt{{{R_{max}(\tau_i)}\over 2}}}} = \sqrt{ {{2 \pi^2} \over {R_{max}(\tau_i)} }},
\]

we have
\[
\Delta (\xi,\tau_i,X,t_i) \leq 9  {{\left(\sqrt{ {{2 \pi^2} \over
{R_{max}(\tau_i)} } }\right)^2\over {t_i-\tau_i}}} = 36 \pi^2.
\]

When $t_i,\tau_i$ large enough, we have
\[
  {  {e^{t_i}-1} \over {e^{\tau_i}-1}} < 2.
\]

Then the Harnack inequality gives
\[
                 R(\xi,\tau_i) \leq  2 e^{9 \pi^2} R(X,t_i)
\]
or
\[
   R(X,t_i) \geq {1 \over 2} e^{-9 \pi^2}\;R_{max}(\tau_i)
\]
for all $X$ in a ball around $\xi$ of radius
\[
             \rho = {{\pi}\over{\sqrt{{{R_{max(\tau_i)}}\over 2}}}}.
\]

By Lemma 7.1, the injectivity radius of the evolved metric at time
$t_i$ is:
 \[ inj_M(t_i) \geq
{{\beta \pi}\over{\sqrt{{{R_{max(t_i)}}\over {2n(n+1)}}}}} \geq
{\beta \rho}\cdot \sqrt{ {n (n+1)} \over 2} .
 \]
 Set $\rho_1 = \rho \cdot\sqrt{ {n (n+1)} \over 2}.\;$  In complex
dimension 2, if
\[
  \displaystyle\;\min_{X \in B_{\beta\rho_1}(\xi)}\;R(X,t_i) > 2 \;r,\]
then
 \[
R(X, t_i) - r  > {1\over 2 } \; R(X, t_i), \qquad \forall\; X \in
B_{\beta\rho_1}(\xi).
 \]

 Therefore, we have (at time $t_i$)
\[
\displaystyle \int_{B_{\beta\rho_1}(\xi)}\; (R(X,t_i) - r)^2
\;{\omega_{\varphi(t_i)}}^n \geq {1\over 4}
   \displaystyle \int_{B_{\beta\rho_1}(\xi)}\; R(X,t_i)^2 \;{\omega_{\varphi(t_i)}}^n  > C.
\]

The last inequality follows from a volume comparison theorem (cf.
\cite{CE75}). This contradicts with the initial assumption
(\ref{eq:sequence}).  Thus,

\[
  \displaystyle\;\min_{X \in B_{\beta\rho_1}(\xi)}\;R(X,t_i) <  2 r.\]

Again, by the Harnack inequality, we have
\[
\begin{array}{lcl}
 \displaystyle\;\max _{B_{\beta\rho}(\xi)}\;  R(X,t_i) & \leq & 2 R_{max}(\tau_i)\\
& \leq  & 4 e^{9 \pi^2}
\displaystyle\;\min_{B_{\beta\rho}(\xi)}\;R(t_i) \\
& < & 8 e^{9 \pi^2} r.
\end{array}
\]

Then the scalar curvature must be uniformly bounded above for this
sequence $t_i \rightarrow \infty$.
\end{proof}

Now combine Lemmas 9.1, 9.2 and 9.3, we can prove the following
theorem:

\begin{theo} In dimension 2 (or Lemma 9.3 holds), then
$R_{max}(t)$ is bounded from above  uniformly along the K\"ahler
Ricci flow.
\end{theo}
\begin{proof} Choose $T = {1\over {16 r}} e^{-9\pi^2}  $ in Lemma 9.2.
Let $\{s_i\},\{t_i\}$ be two sequences as in Lemma 9.2. Then,
Lemma 9.3 implies that $R_{max}(t_i)$ is uniformly bounded by a
constant $8 e^{9\pi^2} r.\;$ Since $s_i \leq t_i + T,$  Lemma 9.1
implies that
 $R_{max}(s_i)$
is uniformly bounded by a constant $ 16 e^{9\pi^2} r.\;$ Since
$s_i$ is an arbitrary sequence of time, the maximal scalar
curvature must be bounded from above uniformly.
\end{proof}

\subsection{ Convergence to K\"ahler-Einstein metrics by sequence}
In this subsection, we want to show that for any integer $l> 0,$
the K\"ahler Ricci flow converges to  a K\"ahler-Einstein metric
in any $C^l$ norm.
 Note that the
limit K\"ahler-Einstein metric may be different when extracting
from a different sequences. We will defer to the next section to
prove that the limit metric is in fact unique.  Let us first
recall
 a theorem by  W. X. Shi \cite{Shi89} (which we have restated in our
 setting):
\begin{theo} \cite{Shi89} Let $(M,g_0)$ be a K\"ahler metric in $M^n$ with bounded sectional
curvature satisfying:
\[  |R_{i\overline{j}k \overline{l}} |^2 \leq k_0,\;\qquad \forall\;i,j,k,l =1,2,\cdots, n.
\]
Then there exists a constant $T(n,k_0)$ which depends only on $n$
and $k_0$ such that the evolution equation
\[\begin{array}{lcl}
  {{\partial g_{i\overline{j}}}\over {\partial t}} & = & g_{i\overline{j}} - R_{i\overline{j}}\qquad {\rm on}\; M,
  \\
   g_{i\overline{j}}(x,0) &  = & {g_0}_{i\overline{j}}(x),\qquad \forall\;x \in M
\end{array}
\]
has a smooth solution in $0 \leq t \leq T(n,k_0),$ and satisfies
the following estimates: For any integer $m\geq 0,$ there exists
constants $c_m >0$ depending only on $n,m,$ and $k_0$ such that
\[
   \displaystyle \sup_{x\in M}\; |\nabla^m R_{i\overline{j}k\overline{l}} | \leq {{c_m}\over{t^m}}, \qquad \forall\; 0\leq t\leq T(n,k_0).
\]

In particular, there exists a constant $c$ such that
\[
   {1\over c} \;g_{i\overline{j}}(x) \leq \tilde{g}_{i\overline{j}}(x) \leq c \;g_{i\overline{j}}(x)
\]
where $\tilde{g}_{i\overline{j}}(x) = g_{i\overline{j}}(x,T).$
\end{theo}
 Combining Shi's theorem with Theorem 9.4, we arrive at the
 following
\begin{theo} The folllowing statements hold along the K\"ahler Ricci flow:

\begin{enumerate}
\item  The injectivity radius has a uniform positive lower bound, and the diameter has
a uniform upper bound.
\item The bisectional curvature and all its derivatives are uniformly
bounded from above over the K\"ahler Ricci flow. In particular,
 the scalar curvature has a uniform upper bound and positive lower bound.
\item $\displaystyle \lim_{t\rightarrow \infty} (R-r) =0.\;$
\item For any integer $l > 0,$ and for any time sequence $t_i
\rightarrow \infty, $ there exists a subsequence of $\{t_i\}$
(still using the same notation) such that the evolved K\"ahler
metrics   converges to a K\"ahler-Einstein metric with constant
bisectional curvature in $C^{l} $ norm.
\end{enumerate}
\end{theo}

\begin{proof} By Theorem 9.4, the bisectioinal curvature $R$ is
uniformly bounded.  Lemma 7.1 then implies that the injectivity
radius has a uniform positive lower bound, which in turns implies
that the Sobolev constant has a uniform upper bound. Since  the
volume is fixed along
  the K\"ahler Ricci flow, the diameters are  bounded uniformly from above.
Repeatedly applying the theorem of Shi, we can show  that all the
derivatives of the sectional curvatures are uniformly bounded
over the entire flow.   In particular, any sequence of metrics
over time must have a subsequence which converges to a limit
metric in $C^{l} $ for any fixed integer $l.\;$ \\

Next we want to show that $\displaystyle \lim_{t \rightarrow
\infty} (R(t) - r) = 0.\;$  We just need to show this for an
arbitrary sequence $s_i \rightarrow \infty.\; $ Lemma 9.2 implies
that there exists another sequence of time $t_i \rightarrow
\infty$ such that
\[
 \displaystyle \lim_{i\rightarrow \infty} \displaystyle \int_M (R(\omega_{\varphi(t_i)}-r)^2
 {\omega_{\varphi(t_i)}}^2=0,\qquad {\rm where}\; t_i \leq s_i
 \leq t_i + T.
 \]

 Combining this with the earlier result of convergence in any $C^l$
 norm, we arrive at
 \[
 \displaystyle \lim_{i \rightarrow \infty}
 (R(\omega_{\varphi(t_i)}-r)= 0.
\]

Thus $\omega_{\varphi(t_i)}$ converges to a K\"ahler-Einstein
metric as $t_i \rightarrow \infty.\;$ Note that all of the $l-$th
derivatives of the evolved metrics are controlled for any integer
$l\geq 0.\;$ Consider the sequence of the K\"ahler Ricci flow from
$t_i$ to $t_i + T.\; $ This sequence of Ricci flow with fixed
length $T$  converges strongly to the K\"ahler Ricci flow of limit
metrics. Since the limit of $\omega_{\varphi(t_i)}$ is a
K\"ahler-Einstein metric, the limiting Ricci flow must be trivial
and all of the limits of sequences of flow from $t_i$ to $t_i +T$
are K\"ahler Einstein metrics. In particular, since $t_i \leq s_i
\leq t_i + T,\;$ we show that $ \displaystyle \lim_{i \rightarrow
\infty} (R(\omega_{\varphi(s_i)}-r)= 0.\;$ Since $\{s_i\}$ is a
sequence chosen randomly,  we then have
\[
 \displaystyle \lim_{t \rightarrow \infty}
 (R(\omega_{\varphi(t)}-r)= 0.
\]
In other words, the limit metric of any sequence  along the
K\"ahler Ricci flow must be of constant scalar curvature.
Consequently, the limit metric of any sequence must be  a
K\"ahler-Einstein metric. Moreover, in $\CC P^n,$ this in turn
implies that the limit metric has  constant bisectional
curvature. In summary,  we have
\[
\displaystyle \lim_{t\rightarrow \infty}\;(R_{i\overline{j}}
-{1\over n} R g_{i\overline{j}}) =0,
\]
and
\[
\displaystyle \lim_{t\rightarrow \infty}\;(R_{i\overline{j}k
\overline{l}} -{1\over {n(n+1)}} R (g_{i\overline{j}}
g_{k\overline{l}} + g_{i\overline{l}} g_{k\overline{j}} )) = 0.
\]
\end{proof}

\section{Exponential convergence}
In the previous section, we prove that the K\"ahler Ricci flow
converges to K\"ahler-Einstein metrics by sequences. Although
limit metrics (from different time sequences) might be isometric
to each other, but certainly not necessarily unique. We want to
show that the limit is unique and the K\"ahler Ricci flow
converges exponentially to this metric. In the 1st subsection,
 we explain how to initialize the K\"ahler potential at time $t=0$ in order to have
 convergence on the K\"ahler potential level. In the second subsection, we prove that
 K\"ahler Ricci flow converges exponentially
 fast to a unique K\"ahler-Einstein metric.

 \subsection{Normalization of initial value}

 Consider the Ricci flow on the K\"ahler potential level,
\begin{equation}
   {{\partial \varphi} \over {\partial t }} =  \log {{\omega_{\varphi}}^n \over {\omega}^n } + \varphi - h_{\omega}.
\label{eq:flowpotential3}
\end{equation}

Define
\[
  c(t) = {1 \over V}\displaystyle \int_M  {{\partial \varphi} \over {\partial t }}
 {\omega_{\varphi}}^n .
\]
We have the following lemma
\begin{lem} Set the inital value of
$\varphi$ at time $0$ so that
\[
c(0) = {1 \over V} \displaystyle \int_0^{\infty} \;
e^{-t}\;\displaystyle \int_M |\nabla {{\partial \varphi} \over
{\partial t }}|^2_{\varphi} {\omega_{\varphi}}^n d\,t < C.
\]
This normalization is appropriate when the K energy has a uniform
lower bound along the K\"ahler Ricci flow. Then, $c(t) > 0\;$ for
all time $t.\;$ We have
\[
\displaystyle \lim_{t\rightarrow \infty}\; c(t) =
 \displaystyle \lim_{t\rightarrow \infty}\; {1 \over V}\displaystyle \int_M  {{\partial \varphi} \over {\partial t }}
 {\omega_{\varphi}}^n = 0.
\]
\end{lem}
\begin{proof}
 A simple calculation yields
\[
  c'(t) = c(t) - {1 \over V} \displaystyle \int_M  |\nabla {{\partial \varphi} \over {\partial t }}|^2_{\varphi}
   {\omega_{\varphi}}^n.
\]

 Define
\[
  \epsilon(t) =  {1 \over V} \displaystyle \int_M  |\nabla {{\partial \varphi} \over {\partial t }}|^2_{\varphi}
   {\omega_{\varphi}}^n.
\]

Since the K energy has a lower bound along the K\"ahler Ricci
flow, we have
\[
  \displaystyle \int_0^{\infty} \epsilon(t) d\,t = {1 \over V} \displaystyle \int_0^{\infty} \displaystyle
   \int_M  |\nabla {{\partial \varphi} \over {\partial t }}|^2_{\varphi} \; {\omega_{\varphi}}^n d\,t < C
\]
for some constant $C.\;$ Now, we normalize our initial value of
$c(t)$ as
\[
\begin{array}{lcl}
  c(0) & = &  {1 \over V} \displaystyle \int_0^{\infty} \epsilon(t) e^{-t} d\,t \\
       & = &    {1 \over V} \displaystyle \int_0^{\infty} \displaystyle \int_M
        |\nabla {{\partial \varphi} \over {\partial t }}|^2_{\varphi}  {\omega_{\varphi}}^n   e^{-t} d\,t \\
        &  \leq & {1 \over V} \displaystyle \int_0^{\infty} \displaystyle \int_M
        |\nabla {{\partial \varphi} \over {\partial t }}|^2_{\varphi}  {\omega_{\varphi}}^n d\,t \\
        &  = & \displaystyle \int_0^{\infty}\; {{d\,\nu}\over {d\,t}} \;d\,t = \nu(0) - \nu(\infty)< C.
\end{array}
\]

This shows that our initial setting is correct. From the equation
for $c(t)$, we have
\[
   (e^{-t} c(t))' = - \epsilon(t) e^{-t}.
\]

Integrating this equation from $0$ to $t,$ we have
\[
\begin{array}{lcl}
  e^{-t} c(t) & = & c(0) -  \displaystyle \int_0^{t} \epsilon(\tau) e^{-\tau} d\,\tau \\ & = & \displaystyle \int_0^{\infty} \epsilon(t) e^{-t} d\,t -  \displaystyle \int_0^{t} \epsilon(\tau) e^{-\tau} d\,\tau \\
  & = &  \displaystyle \int_t^{\infty} \epsilon(\tau) e^{-\tau} d\,\tau.
\end{array}
\]

Thus
\begin{eqnarray}
   c(t) & = & e^t  \displaystyle \int_t^{\infty} \epsilon(\tau) e^{-\tau} d\,\tau \nonumber  \\
   & = &  \displaystyle \int_t^{\infty} \epsilon(\tau) e^{-(\tau-t)}
   d\,\tau \nonumber
  \\ & \leq &  \displaystyle \int_t^{\infty} \epsilon(\tau)  d\,\tau \rightarrow 0.
  \label{eq:expo}
\end{eqnarray}

Note that $c(t) > 0$ for all time.  In conclusion, we have
\begin{equation}
 \displaystyle \lim_{t \rightarrow \infty} c(t) = \displaystyle \lim_{t \rightarrow \infty}
  \displaystyle \int_M  {{\partial \varphi} \over {\partial t }}  {\omega_{\varphi}}^n  = 0.
\label{eq:initialnormalization}
\end{equation}
\end{proof}
\subsection{Exponential convergence}
In this subsection, we assume that the evolved K\"ahler metrics
$\omega_{\varphi(t)}$ converge to a K\"ahler-Einstein metric in at
least $C^3-$norm. We then show that the flow must converge to a
unique K\"ahler-Einstein metric exponentially fast.

Recall that the K\"ahler Ricci flow equation:
\[
   {{\partial \varphi} \over {\partial t}} = \ln {{\omega_{\varphi}^n} \over {\omega^n}}
                                + \varphi -h_{\omega}.
\]

Since the evolved K\"ahler metrics converge to K\"ahler-Einstein
metrics by sequences in any $C^k$ norm, then we have the following

\begin{enumerate}
\item Modulo constants, we have
$$ \displaystyle \lim_{t \rightarrow \infty} {{\partial \varphi} \over {\partial t}} =
 \displaystyle \lim_{t \rightarrow \infty} (\ln {{\omega_{\varphi}^n} \over {\omega^n}}
                                + \varphi -h_{\omega}) =0. $$

 This together with the normalization
of the initial value (see Lemma 10.1), we have
\[
\displaystyle \lim_{t \rightarrow \infty} {{\partial \varphi}
\over {\partial t}} =
 \displaystyle \lim_{t \rightarrow \infty} (\ln {{\omega_{\varphi}^n} \over {\omega^n}}
                                + \varphi -h_{\omega}) =0.
\]
\item The eigenspace of $\omega_{\varphi(t)}$ converges to the eigenspace of a K\"ahler-Einstein
metric. Notice that in a fixed K\"ahler class, all
K\"ahler-Einstein metrics are isometric to each other so that
they have the same spectrum.

\item The eigenvalues of $\omega_{\varphi(t)}$ converge to  the eigenvalues
of some K\"ahler Einstein metrics. Note that the second eigenvalue
of a K\"ahler Einstein metric is strictly bigger than $1.\;$

\end{enumerate}

\begin{prop} There exists
 a positive number $\alpha > 0$ and constant $C > 0 $ such that
\[
    \displaystyle \int_M \left({{\partial \varphi} \over {\partial t}} - c(t)\right)^2
                     \omega_{\varphi}^n \leq C \;e^{- \alpha\; t}.
\]
Moreover, for every integer $l > 0,$ there exists a constant
$C_l$ such that
\[
 \displaystyle \int_M \mid D^l\;\left({{\partial
\varphi(t)} \over {\partial t}} -
   c(t)\right)\mid^2_{\varphi}\;\; \omega_{\varphi(t)}^n
   \leq C_l  \; e^{-\alpha t}.
\]
\end{prop}
First we want to prove a corollary of this proposition
\begin{cor} There exists a uniform constant $C$ such that
\[
  0 < c(t) \leq C\; e^{-\alpha \;t}, \qquad \forall\;t > 0.
\]
\end{cor}
\begin{proof}
Recall that

 \[ \epsilon(t) = {1\over V} \displaystyle \int_M \mid
\p \;\left({{\partial \varphi(t)} \over {\partial t}} -
   c(t)\right)\mid^2_{\varphi}\;\; \omega_{\varphi(t)}^n \leq C_1
   \; e^{-\alpha\; t},
 \]
 where $C_1$ is some uniform constant.
 Plugging this into Formula (\ref{eq:expo}), we obtain
 \[
 c(t) =\displaystyle \int_t^{\infty} \epsilon(\tau) e^{-(\tau-t)} d\,\tau
      \leq {C_1 \over {1 + \alpha}}  \; e^{t} e^{-(1+\alpha)t} = {C_1 \over {1 + \alpha}} \;e^{-\alpha\;t}.
 \]
\end{proof}
Next we return to prove Proposition 10.2.

\begin{proof} Differentiating the K\"ahler Ricci flow with respect to time $t,$
\begin{equation}
  {{\partial^2 \varphi} \over {\partial t^2}} = \Delta_{\varphi} {{\partial \varphi} \over {\partial t}}
           + {{\partial \varphi} \over {\partial t}}.
           \label{eq:expodecay0}
\end{equation}

Put $\mu(t) =  \displaystyle \int_M \left({{\partial \varphi}
\over {\partial t}} - c(t)\right)^2
                     \omega_{\varphi}^n.\;$ Then
\[\begin{array}{lcl}
   & & {{d\, \mu(t)}\over {d\,t}}\\
   & = & 2 \displaystyle \int_M \left({{\partial \varphi} \over {\partial t}} - c(t)\right)
   \left({{\partial^2 \varphi} \over {\partial t^2}} - c'(t)\right)
                     \omega_{\varphi}^n  + \displaystyle \int_M \left({{\partial \varphi} \over {\partial t}} - c(t) \right)^2
                     \Delta_{\varphi}  {{\partial \varphi} \over {\partial t}}
                      \omega_{\varphi}^n \\
                      & = & 2 \displaystyle \int_M \left({{\partial \varphi} \over {\partial t}} - c(t)\right) \left( \Delta_{\varphi} {{\partial \varphi} \over {\partial t}}
           + {{\partial \varphi} \over {\partial t}} - c(t)\right)
                     \omega_{\varphi}^n  + \displaystyle \int_M \left({{\partial \varphi} \over {\partial t}} - c(t)\right)^2
                     \Delta_{\varphi}  {{\partial \varphi} \over {\partial t}}
                      \omega_{\varphi}^n\\
    & = & -2 \displaystyle \int_M (1 + {{\partial \varphi} \over {\partial t}} - c(t))\mid \nabla \left({{\partial \varphi} \over
     {\partial t}} - c(t)\right)\mid^2_{\varphi} \omega_{\varphi}^n
       + 2 \displaystyle \int_M \left({{\partial \varphi} \over {\partial t}} - c(t)\right)^2
                     \omega_{\varphi}^n.
                    \end{array} \]
Here we have used the fact $\displaystyle \int_M \left({{\partial
\varphi} \over {\partial t}} - c(t)\right)
 \omega_{\varphi}^n = 0$ twice. Since $ \displaystyle \lim_{t \rightarrow \infty} {{\partial \varphi} \over {\partial t}} =
  \displaystyle \lim_{t \rightarrow \infty} c(t) = 0$ for any $\epsilon > 0,$ and for $t$
  large enough, we have
  \begin{eqnarray}
  & & {{d\, \mu(t)}\over {d\,t}}  \nonumber\\
  & \leq & -  2 (1 -\epsilon) \displaystyle \int_M \mid \nabla \left({{\partial \varphi} \over {\partial t}} - c(t)\right)\mid^2_{\varphi} \omega_{\varphi}^n
    \nonumber    \\
       & & \qquad \qquad+ 2 \displaystyle \int_M \left({{\partial \varphi} \over {\partial t}} - c(t)\right)^2
                     \omega_{\varphi}^n.
                     \label{eq:decay0}
  \end{eqnarray}

If the first eigenvalue of $\omega_{\varphi(t)}$ converges to
$1,$ it appears to be quite difficult from the previous
inequality to derive any control on  ${{d\, \mu(t)}\over
{d\,t}}.\;$
  Denote the first, second eigenvalue of a K\"ahler Einstein metric as
  $\lambda_1 < \lambda_2.\;$ Lemma 6.3 implies that $\lambda_1
  \geq 1$ and the equality holds if and only if the space of holomorphic vector
  fields $\e (M) $ is non-trivial. In case of $\e (M) = 0,$ we
  have $\lambda_1 > 1.\;$ For $t$ large enough, all eigenvalues of
  $\omega_{\varphi(t)}$ will be bigger than ${{\lambda_1 + 1}\over 2} >
  1.\;$  Therefore,
 \[
   \displaystyle \int_M \mid \nabla \left({{\partial \varphi} \over {\partial t}} - c(t)\right)
   \mid^2_{\varphi} \omega_{\varphi}^n   \geq   {{\lambda_1 + 1} \over 2} \displaystyle \int_M
   \; \left({{\partial \varphi} \over {\partial t}} - c(t)\right)^2
   \omega_{\varphi}^n.
 \]

 Plugging  this into inequality (\ref{eq:decay0}), we obtain
   \[
   \begin{array}{lcl}
   {{d\, \mu(t)}\over {d\,t}} & \leq &   -  2  (1-\epsilon) {{\lambda_1 + 1} \over 2}  \displaystyle
      \int_M \left({{\partial \varphi} \over {\partial t}} -c(t)\right)^2  \omega_{\varphi}^n
       + 2 \displaystyle \int_M \left({{\partial \varphi} \over {\partial t}} - c(t)\right)^2
                     \omega_{\varphi}^n \\
   & \leq & - \alpha \displaystyle \int_M \left({{\partial \varphi} \over {\partial t}} - c(t)\right)^2
                     \omega_{\varphi}^n = - \alpha \;\mu(t),
   \end{array}
   \]
   where $\alpha = 2 (1 -\epsilon) {{\lambda_1 + 1} \over 2} - 2. \;$ Choose $\epsilon > 0$
   to be small enough, we have $\alpha > 0.\;$  It is straightforward to prove that
    there exists a uniform constant $C$ such that
   \[
     \mu(t) = \displaystyle \int_M \left({{\partial \varphi} \over {\partial t}} - c(t)\right)^2
                     \omega_{\varphi}^n  \leq C \; e^{-\alpha t}.
   \]

 On the other hand, if $\e (M) \neq 0, $ then $\lambda_1 = 1$ and
 the first eigenvalue of $\omega_{\varphi(t)}$ converges to $1.\;$
 The inequality (\ref{eq:decay0}) gives us little control of the growth of $\mu(t)$.
 However, the Futaki invariant comes to rescue: Let $X$ be any holomorphic vector
  field, then (in a K\"ahler-Einstein manifold)
  \[\begin{array}{lcl}
  0 & = &  f_M(X,\omega_{\varphi}) \\
    & =  & \displaystyle \int_M X\left(\ln {{\omega_{\varphi}^n} \over {\omega^n}}
                                + \varphi -h_{\omega}\right) \omega_{\varphi}^n \\
    & = &  \displaystyle \int_M X\left({{\partial \varphi} \over {\partial t}} -c(t)\right)
    \omega_{\varphi}^n  = - \displaystyle \int_M \Delta_{\varphi} \theta_X \cdot
    \left({{\partial \varphi} \over {\partial t}} -c(t)\right)
    \omega_{\varphi}^n,
  \end{array}
  \]
  where $L_X (\omega_{\varphi}) = \sqrt{-1} \partial \overline{\partial} \theta_X$ as defined in Section 6.
  If $\omega_{\varphi}$ were already a K\"ahler-Einstein metric, then the above inequality would imply that
  \[
    \displaystyle \int_M \Delta_{\varphi} (\theta_X) \cdot
    \left({{\partial \varphi} \over {\partial t}} - c(t)\right)  \omega_{\varphi}^n
    = - \displaystyle \int_M \theta_X \cdot
    \left({{\partial \varphi} \over {\partial t}} -c(t)\right)  \omega_{\varphi}^n =
    0.
  \]

  This in its turn would imply that $\left({{\partial \varphi} \over {\partial t}} -c(t)\right)$ is perpendicular to
  the first eigenspace of $\omega_{\varphi}.\;$ And that would  give us the desired estimate from the
  inequality (\ref{eq:decay0}).   Unfortunately, $\omega_{\varphi}$ is not
  a K\"ahler-Einstein metric. However, $\omega_{\varphi(t)}$ is at least $C^3$ close
  to a K\"ahler-Einstein metric as $t \rightarrow \infty;\;$ and this shall be sufficient to
  derive the exponential convergence.  Note that the eigenvalues of $\omega_{\varphi(t)}$ shall converges
  to the eigenvalues of  K\"ahler Einstein metrics. For any fixed $\epsilon >
  0,$ and for $t$ large enough, the eigenvalues of $\omega_{\varphi(t)}$ must be either
    in $(1-\epsilon, 1+ \epsilon)$ or are strictly bigger than ${{\lambda_2 + 1}\over 2} > 1 + \epsilon.\;$
    Denote the set of all eigenspaces of $\omega_{\varphi(t)}$  whose eigenvalues
    are between $(1-\epsilon,1+\epsilon) $ as $\Lambda_{{\rm small}}(\omega_{\varphi}).\;$
    Then $\Lambda_{{\rm small}}(\omega_{\varphi}) $  converges to the first eigenspace of
    some K\"ahler Einstein metrics. Moreover, by Lemma 6.3, $\{\triangle_{\varphi(t)} \theta_X\;\mid \; X \;\in
    \;\e ({\rm M}) \}$ converges to the first eigenspace of the limit K\"ahler-Einstein metric
    space. Thus, $\{\triangle_{\varphi(t)} \theta_X\;\mid \;  X \;\in
    \;\e ({\rm M}) \}$ is essentially $\Lambda_{\rm small}(\omega_{\varphi(t)}),$ where possible
     error terms become as small as needed when $t \rightarrow \infty.\;$ In
other words, the vanishing of Futaki
     invariant implies that the
   projection of ${{\partial \varphi} \over {\partial t}} -c(t)$
   into
the  eigenspace $\Lambda_{{\rm small}}(\omega_{\varphi}) $ is
very small (compare to the size of
   $ {{\partial \varphi} \over {\partial t}} -c(t)$).  Namely, we have
   \[
     {{\partial \varphi} \over {\partial t}} -c(t) = \varrho +
     \varrho^{\perp},
   \]
   where $\varrho \in \Lambda_{\rm small} (\omega_{\varphi})$
   and $\varrho^{\perp}\perp \Lambda_{\rm small} (\omega_{\varphi}).\;$ For $t$ large enough, we have
   \[
     \displaystyle \int_M \varrho^2  \omega_{\varphi}^n \leq \epsilon \displaystyle
      \int_M \left({{\partial \varphi} \over {\partial t}} -c(t)\right)^2  \omega_{\varphi}^n,
   \]
   and
   \[
   \displaystyle \int_M {\varrho^\perp}^2  \omega_{\varphi}^n \geq (1 - \epsilon) \displaystyle
      \int_M \left({{\partial \varphi} \over {\partial t}} -c(t)\right)^2  \omega_{\varphi}^n.
   \]

   Notice that the eigenvalue of $\omega_{\varphi(t)}$ corresponds to ${\Lambda_{\rm small} (\omega_{\varphi(t)})}^{\perp}$
     are always bigger than ${{\lambda_2 + 1} \over 2} > 1\;$ when
     $t$ large enough.
   Therefore,
   \[
   \begin{array}{lcl}
   \displaystyle \int_M \mid \nabla \left({{\partial \varphi} \over {\partial t}} - c(t)\right)
   \mid^2 \omega_{\varphi}^n  & \geq  & {{\lambda_2 + 1} \over 2} \displaystyle \int_M  {\varrho^\perp}^2
   \omega_{\varphi}^n \\
   & \geq & (1-\epsilon) {{\lambda_2 + 1} \over 2}  \displaystyle
      \int_M \left({{\partial \varphi} \over {\partial t}} -c(t)\right)^2  \omega_{\varphi}^n.
      \end{array}
   \]

   Plugging  this into inequality (\ref{eq:decay0}), we obtain
   \[
   \begin{array}{lcl}
   {{d\, \mu(t)}\over {d\,t}} & \leq &   -  2 (1 -\epsilon) (1-\epsilon) {{\lambda_2 + 1} \over 2}  \displaystyle
      \int_M \left({{\partial \varphi} \over {\partial t}} -c(t)\right)^2  \omega_{\varphi}^n
       + 2 \displaystyle \int_M \left({{\partial \varphi} \over {\partial t}} - c(t)\right)^2
                     \omega_{\varphi}^n \\
   & \leq & - \alpha \displaystyle \int_M \left({{\partial \varphi} \over {\partial t}} - c(t)\right)^2
                     \omega_{\varphi}^n = - \alpha \;\mu(t),
   \end{array}
   \]
   where $\alpha = 2 (1 -\epsilon)^2 {{\lambda_2 + 1} \over 2} - 2. \;$ Again, we choose $\epsilon > 0$
   to be small enough, we have $\alpha > 0.\;$  It is straightforward to prove that
    there exists a uniform constant $C$ such that
   \[
     \mu(t) = \displaystyle \int_M \left({{\partial \varphi} \over {\partial t}} - c(t)\right)^2
                     \omega_{\varphi}^n  \leq C \; e^{-\alpha t}.
   \]

This proves the first part of Proposition 10.2. Next we want to
prove the exponential convergence for all derivatives.  For any
integer $l>0,$ consider the $L^2$ norm of $l-$th derivatives
($l\geq 1$) of
 $\left({{\partial \varphi} \over {\partial t}} - c(t)\right): $
   \[
   \begin{array}{lcl}
   \mu_l(t) & = &\displaystyle \int_M \mid D^l\;\left({{\partial \varphi(t)} \over {\partial t}} -
   c(t) \right)\mid^2_{{\varphi(t)}}\;\; \omega_{\varphi(t)}^n \\
   & = & \displaystyle \int_M \mid D^l\;
   {{\partial \varphi(t)} \over {\partial t}}
  \mid^2_{{\varphi(t)}}\;\; \omega_{\varphi(t)}^n.
\end{array}
   \]
   Since $\displaystyle \lim_{t\rightarrow \infty} \left({{\partial \varphi} \over {\partial t}} - c(t) \right)
   = 0$ and since the K\"ahler Ricci flow converges to some K\"ahler-Einstein metrics
    in any $C^k$ norm (Theorem 9.6) for any integer $k>0$, we have
   \[
    \displaystyle \lim_{t \rightarrow \infty} \mu_l(t) = 0.
   \]
   Finally, we want to show that it is exponentially decay along the
   K\"ahler Ricci flow.
   Using the equation (\ref{eq:expodecay0}), and the fact that all
   the derivatives of curvature are uniformly bounded, we have
   ($l\geq 1$)
   \[
   \begin{array}{lcl} & &{{d \mu_l(t) }\over {d\,t}}\\ & = & \displaystyle \int_M {{\partial }\over {\partial t}}
   \;\mid D^l\;{{\partial \varphi(t)} \over {\partial t}} \mid^2_{\varphi}
   \omega_{\varphi(t)}^n  + \displaystyle \int_M \mid D^l\;{{\partial \varphi(t)} \over {\partial t}}
   \mid^2_{\varphi}\;\;\triangle_{\varphi(t)}\;{{\partial \varphi(t)} \over {\partial t}}\;
    \omega_{\varphi(t)}^n \\
    & \leq  & - 2\;  \displaystyle \int_M \mid D^{l+1}\;{{\partial \varphi(t)} \over {\partial t}}
    \mid^2_{\varphi}\;\; \omega_{\varphi(t)}^n  +
c(n,l) \displaystyle \int_M \mid D^l\;{{\partial
\varphi(t)} \over {\partial t}} \mid^2_{\varphi}\;\;\omega_{\varphi(t)}^n.\\
   & \leq & - 2\;  \displaystyle \int_M \mid D^{l+1}\;{{\partial \varphi(t)} \over {\partial t}} \mid^2_{\varphi}\;
   \; \omega_{\varphi(t)}^n
   \\
   & &  +
   c(n,l) \left( \epsilon   \displaystyle \int_M \mid D^{l+1}\;{{\partial \varphi(t)} \over {\partial t}} \mid^2_{\varphi}\;\; \omega_{\varphi(t)}^n +
c(\epsilon)   \displaystyle \int_M \left({{\partial \varphi(t)}
\over {\partial t}} - c(t)\right)^2 \; \omega_{\varphi(t)}^n  \right) \\
   & = & - (2-c(n,l) \epsilon) \;  \displaystyle \int_M \mid D^{l+1}\;{{\partial \varphi(t)} \over {\partial t}}
   \mid^2_{\varphi}\;\; \omega_{\varphi(t)}^n \\
   & & \qquad \qquad
   + c(n,l) c(\epsilon) \displaystyle \int_M \left({{\partial \varphi(t)} \over {\partial t}} -
   c(t)\right)^2\; \omega_{\varphi(t)}^n.
   \end{array}
  \]
  In the first inequality, we use integration by parts and the
  fact that all of the $l-$th derivatives of the metrics are
  uniformly bounded. In the second to the last inequality, we have
  used an interpolation formula where $C(\epsilon)$ is the interpolation
  constant. Choose $\epsilon$ to be small enough so that
  \[
  2-c(n,l) \epsilon > 0.
  \]
  Then, we have
  \[
{{d \mu_l(t) }\over {d\,t}} \leq c(n,l) c(\epsilon) \displaystyle
\int_M \left({{\partial \varphi(t)} \over {\partial t}} -
   c(t)\right)^2\; \omega_{\varphi(t)}^n \leq C_l \;e^{-\alpha t}.
  \]
Here $C_l = c(n,l) c(\epsilon) C < \infty.\; $ Integrating the
above inequality from $t$ to $\infty,$ we arrive at the desired
estimates :
\[
\mu_l(t) = \displaystyle \int_M \mid D^l\;\left({{\partial
\varphi(t)} \over {\partial t}} -
   c(t)\right)\mid^2_{\varphi}\;\; \omega_{\varphi(t)}^n
   \leq C_l  \; e^{-\alpha t},  \]
where we have used the fact that $\displaystyle
\lim_{t\rightarrow \infty} \mu_l(t) = 0.\;$

\end{proof}

 Note that $\omega_{\varphi(t)}$ have uniform
positive lower bound on injective radius and a uniform positive
bound for Sobolev constant. Combining the above inequality  and
Sobolev embedding theorem, we arrive at
\begin{equation}
\mid D^l\;\left({{\partial \varphi(t)} \over {\partial t}} -
   c(t)\right)\mid^2_{\varphi} \leq c_{l}
   \; e^{-\alpha t},
   \label{ex:last}
\end{equation}
where $c_l$ is another set of uniform constants. In particular
when $l=0,$ we have

\[
\mid {{\partial \varphi(t)} \over {\partial t}} -
   c(t) \mid < c_0 \; e^{-\alpha t}.
\]
Combining this inequality with Corollary 10.3, we have
\[
\mid {{\partial \varphi(t)} \over {\partial t}} \mid < c_0 \;
e^{-\alpha t}
\]
where $c_0$ might be some new constant. Thus there exists a
unique K\"ahler potential $\varphi(\infty)$ such that
\[
\mid \varphi(t) - \varphi(\infty)\mid \leq c_0 \; e^{-\alpha t}.
\]

 From here, we can easily obtain that $\omega_{\varphi(t)} (0 \leq t \leq \infty) $
 are mutually equivalent, i.e., there exists a uniform constant $c>1$ such that
 \[
  {1\over c}   \omega_{\varphi(\infty)} \leq \omega_{\varphi(t)}
  \leq c\; \omega_{\varphi(\infty)}, \qquad \forall \; t \geq 0.
 \]
 Here $\varphi(\infty)$ is the  unique K\"ahler Einstein
 metric (arisen from the limit of the K\"ahler Ricci flow).
 Combining this with inequalities (\ref{ex:last}), we can
   easily imply that

   \[
\mid D^l\;\left({{\partial \varphi(t)} \over {\partial t}} -
   c(t)\right)\mid^2_{\varphi(\infty)} \leq c_{l}
   \; e^{-\alpha t},
   \]

   Thus, $\varphi(t)$ converges exponentially fast to a unique K\"ahler-Einstein metric in ${\cal P}(M,\omega)$
   in any $C^l$ norm. We then prove the following proposition
   \begin{prop} For any integer $l>0,\;{{\partial \varphi} \over {\partial t}} $
   converges exponentially fast to $0\;$ in any $C^l$ norm. Furthermore, the
   K\"ahler Ricci flow converges exponentially fast to a unique
   K\"ahler Einstein metric on any K\"ahler-Einstein surfaces.
   \end{prop}
\section{Concluding Remarks}
Now we prove our main Theorem 1.1 and Corollary  1.2.
\begin{proof} Theorem 1.1 follows from Proposition 10.4. Next we
want to prove Corollary 1.2.  For any K\"ahler metric  in the
canonical K\"ahler class  such that it has non-negative
bisectional curvature on $M$ but positive bisectional curvature
at least at one point, we apply the K\"ahler Ricci flow to this
metric. According to Theorem 2.2,  the bisectional curvature of
the evolved metric  is strictly positive. By our theorem 1.1, the
K\"ahler Ricci flow converges exponentially  to a unique K\"ahler
Einstein metric with constant positive bisectional curvature.
Thus, any K\"ahler metric with nonnegative bisectional curvature
on $M$ and positive at least at one point is path connected to a
K\"ahler-Einstein metric with positive bisectional curvature.
Note that all the K\"ahler-Einstein metrics are path connected by
automorphisms \cite{Ma85}. Therefore, the space of all K\"ahler
metrics with nonnegative bisectional curvature on $M$ and
positive at least at one point, is path connected. Similarly,
using Theorem 2.3 and our Theorem 1.1, we can show that all of the
K\"ahler metrics with nonnegative curvature operator on $M$  and
positive at least at one point is path connected. Note that the
nonnegative curvature operator implies the nonnegative
bisectional curvature.
\end{proof}

\begin{rem} Combining our main theorem 1.1 and Theorem 2.2, 2.3, we
can easily generalize Corollary 1.2 to the case that the
bisectional curvature (or curvature operator) is only assumed to
be non-negative.
\end{rem}
Next we want to propose some future problems.
\begin{q} As our remark 1.4 indicates, what we really need is the
positivity of Ricci curvature along the K\"ahler Ricci flow.
However, it is not expected that the positivity of Ricci
curvature is preserved under the K\"ahler Ricci flow except on
Riemann surfaces. The positivity of bisectional curvature is a
technical assumption to assure the positivity of Ricci curvature.
It is very interesting to extend Theorem 1.1 to metrics without
the assumption on bisectional curvature.
\end{q}

\begin{q} Is the positivity of the sectional curvature preserved under
K\"ahler Ricci flow?
\end{q}

\section{Appendix: Sectional curvature and bisectional curvature}
In this appendix, we want to derive a formula which expresses the
sectional curvature in terms of the bisectional curvature. It is
well known  that in a K\"ahler manifold, these two types of
curvature tensors determine each other uniquely. For the reader's
convenience, we included such a formula here.

We first explain some basic concepts  of the sectional curvature
and the bisectional curvature. Let $u,v,w,x $ be any four tangent
vectors in $M.\;$ Suppose $R(u,v,w,x)$ is the Riemannian curvature
tensor.
 Then
 \[
   R(u,v,Jw,Jx) = R(u,v,w,x)
 \]
 where $J$ is the  complex structure of $M.\;$
 Because of splitting $T_{C} M = T^{1,0}M \oplus T^{0,1}M$ into $\pm \sqrt{-1}$ eigenspaces
of $J$, we can deduce that $R(u,v,w,x) = 0$ unless $w$ and $x$
are of different type. We will use this property strongly in the
this appendix.   Suppose  $x \perp y$ are  two unit tangent
vectors of $M$. Denote the sectional curvature on the plane $x,y$
as $K(x,y).\;$  Set now
\[
  u = {1\over {\sqrt{2}}} ( x - \sqrt{-1} J x),\qquad v = {1\over {\sqrt{2}}} ( y - \sqrt{-1} J
  y).
\]

 If $y \bot J \;x,
$ then
\begin{equation}
R(u,\overline{u},v,\overline{v}) = R(x,y,y,x) + R(x,J\; y, J\;
y,x). \label{eq:section}
\end{equation}

If $y=Jx,$ then
\[R(u,\overline{u},v,\overline{v}) = R(x,J\; x,J\; x,x).
\]

This means that the bisectional curvature and the sectional
curvature are the same on any holomorphic plane. Now, we seek a
formula which expresses the sectional curvature in terms of the
bisectional curvature.

\begin{theo} If $w_1, w_2$  are two mutually  perpdicular real vectors
 in $TM$ such that the two complex
planes   spanned by $w_1$ and $w_2$ respectively are either
perpendicular to each other or are identical, then the sectional
curvature of the real plane spanned by  these two vector fields is
\[
K(w_1,w_2) = {1\over 4}\left(R(A,\overline{A}, A,\overline{A}) -2
R(B, \overline{B}, A, \overline{A}) +
R(B,\overline{B},B,\overline{B}) \right)
\]
where $A = {1\over \sqrt{2}} (u_1 +  u_2)$ and $B = {1\over
\sqrt{2}} (u_1 - u_2)$ and
\[
u_1 = {1\over \sqrt{2}} (w_1 - \sqrt{-1}J w_1),\qquad {\rm and}
\qquad u_2 = {1\over \sqrt{2}} (-J w_2 -\sqrt{-1} w_2).
\]
\end{theo}

An immediate corollary is
\begin{cor} If the bisectional curvature is less than $1,$ then the sectional curvature is
less than $2.\;$
\end{cor}
Let us prove the corollary first.

\begin{proof} If the two complex plane spanned by $w_1$ and $w_2$
are identical, then we must have (if necessary, we can change
$w_2$ to $-w_2$):
\[
    w_2 =  J w_1 \qquad {\rm and} \qquad w_1 = - J w_2.
\]

 Then $u_1 = u_2,\;$ which in turn implies that  $A = 2\; u_1$ and $B = 0.\;$
 Therefore,
 \[
  K(w_1,w_2) = R(u_1,\bar u_1,u_1,\bar u_1) \leq 2.
 \]

 On the other hand, if the two complex planes spanned by $w_1$ and $w_2$
 are mutually perpendicular, then $A, \;B$ are both unit vectors and
 $A \perp B.\;$  Thus \footnote{In this calculation, if the metric has constant bisectional curvature,
 then $R(A,\bar A, B, \bar B) = 1$ and this yields that
 \[
   K(w_1, w_2) = {1\over 4}\left(2-2+ 2 \right) = {1\over 2}.
 \]

 Therefore in case of constant bisectional curvature $1,$ the sectional curvature
 is between ${1\over 2} $ and $2.\;$ }
 \[
 \begin{array}{lcl}
 K(w_1,w_2) & = & {1\over 4}\left(R(A,\overline{A}, A,\overline{A}) -2
R(B, \overline{B}, A, \overline{A}) +
R(B,\overline{B},B,\overline{B}) \right)\\
& \leq & {1\over 4}\left(R(A,\overline{A}, A,\overline{A}) +
R(B,\overline{B},B,\overline{B}) \right) \\
& \leq & {1\over 4}\left(2 + 2 \right) = 1.
 \end{array}\]

In the first inequality in this calculation, we used the fact that
the bisectional curvature is positive:
\[
  R(A,\bar A, B, \bar B) \geq 0.
\]

In conclusion, we prove that the sectional curvature must be less
than 2.
\end{proof}
Now we are ready to give a proof of the main theorem in this
Appendix.
\begin{proof}

 In a local coordinate, let us choose an orthonormal  basis $e_1,e_2,\cdots, e_{2n}$ such that
 $J e_i = e_{n+i}$ for $i=1,2,\cdots, n$ and set
 \[
  e_1 = w_1,  \qquad {\rm and}\qquad e_2 = - J w_2.
 \]

 Let $u_i = {1 \over \sqrt{2}} (e_i - \sqrt{-1} J e_i).\;$
 Then $\{u_i\}$ is a unitary basis.  Conversely, we have
 \[  e_i = {1
\over \sqrt{2}} (u_i + \overline{u_i}),\qquad {\rm  and}\;\;  J
e_i = {{\sqrt{-1}} \over \sqrt{2}} (u_i - \overline{u_i}).\]

 Then,
\begin{eqnarray} 4 K(e_1, J e_i) & = & - 4 \;R(e_1, J e_i,e_1, J e_i)  \nonumber \\
  & = & - 4\; R({1 \over \sqrt{2}} (u_1 + \overline{u_1}),{\sqrt{-1} \over \sqrt{2}} (u_i - \overline{u_i}),
        {1 \over \sqrt{2}} (u_1 + \overline{u_1}),{\sqrt{-1} \over \sqrt{2}} (u_i - \overline{u_i}))
        \nonumber \\
  & = &  \{R(u_1, \overline{u_i}, u_1, \overline{u_i}) - R(u_1, \overline{u_i}, \overline{u_1},u_i)
       - R(\overline{u_1}, u_i, u_1, \overline{u_i} ) + R(\overline{u_1}, u_i,\overline{u_1},u_i) \} \nonumber\\
  & = & \{R(u_1, \overline{u_i}, u_1, \overline{u_i})
        +R(u_1, \overline{u_i}, u_i,\overline{u_1})
        \nonumber \\
       & & \qquad \qquad  \qquad \qquad+ R(u_i,\overline{u_1}, u_1, \overline{u_i} ) + R(u_i,\overline{u_1}, u_i,\overline{u_1})
       \}.  \label{eq:sect1}
\end{eqnarray}

Let $ v $ be any vector in $T^{1,0} M.\;$ For any $\theta = \pm
1,$ we have
\[
  \begin{array}{lcl} R(u_1 + \theta u_i, \overline{u_1} + \theta \overline{u_i}, v, \overline{v}) & = &
    R(u_1, \overline{u_1}, v,\overline{v}) + \theta^2 R(u_i, \overline{u_i}, v,\overline{v})
    \\ & & \qquad + \theta \left( R(u_1, \overline{u_i}, v, \overline{v}) + R(u_i, \overline{u_1}, v, \overline{v})\right).
  \end{array}
\]

Thus,
\begin{eqnarray} & &
2\{R(u_1, \overline{u_i}, v, \overline{v}) + R(u_i,
\overline{u_1}, v, \overline{v})\} \nonumber  \\
 & & \qquad \qquad =  R(u_1 + u_i, \overline{u_1} +
\overline{u_i}, v, \overline{v})
  - R(u_1 -  u_i, \overline{u_1} - \overline{u_i}, v, \overline{v}).
\label{eq:sect2}
\end{eqnarray}

Let $v = u_1 + \varsigma u_i$ and $\varsigma = \pm 1.\;$ Then
\[\begin{array}{lcl}
R(u_1, \overline{u_i}, u_1 + \varsigma u_i, \overline{u_1} +
\varsigma \overline{u_i}) & = &
R(u_1, \overline{u_i}, u_1 , \overline{u_1} ) + \varsigma^2 R(u_1, \overline{u_i}, u_i, \overline{u_i}) \\
& & \qquad +  \varsigma \{ R(u_1, \overline{u_i}, u_1,
\overline{u_i}) + R(u_1, \overline{u_i}, u_i, \overline{u_1})\}.
\end{array}
\]

Thus,
\begin{eqnarray}
R(u_1, \overline{u_i}, u_1 +  u_i, \overline{u_1} +
\overline{u_i}) - R(u_1, \overline{u_i}, u_1 - u_i,
\overline{u_1} -  \overline{u_i})\nonumber \\ \qquad \qquad = 2 \{
R(u_1, \overline{u_i}, u_1, \overline{u_i}) + R(u_1,
\overline{u_i}, u_i, \overline{u_1})\}. \label{eq:sect3}
\end{eqnarray}

Switch $i$ and $1$ in the previous formula, we obtain
\begin{eqnarray}
R(u_i, \overline{u_1}, u_1 +  u_i, \overline{u_1} +
\overline{u_i}) - R(u_i, \overline{u_1}, u_1 - u_i,
\overline{u_1} -  \overline{u_i}) \nonumber \\  \qquad \qquad = 2
\{ R(u_i, \overline{u_1}, u_1, \overline{u_i}) + R(u_i,
\overline{u_1}, u_i, \overline{u_1})\}. \label{eq:sect4}
\end{eqnarray}

 Adding equation (\ref{eq:sect3}) and (\ref{eq:sect4}) together, and using equation (\ref{eq:sect2}), we obtain
\[
\begin{array}{l} 2 \{ R(u_1, \overline{u_i}, u_1, \overline{u_i}) + R(u_1, \overline{u_i}, u_i, \overline{u_1})
       + R(u_i, \overline{u_1}, u_1, \overline{u_i}) + R(u_i, \overline{u_1}, u_i, \overline{u_1})\} \\
= R(u_1, \overline{u_i}, u_1 +  u_i, \overline{u_1} +
\overline{u_i}) - R(u_1, \overline{u_i}, u_1 - u_i,
\overline{u_1} -  \overline{u_i}) \\ \qquad  \qquad \qquad \qquad
+ R(u_i, \overline{u_1}, u_1 +  u_i, \overline{u_1} +
\overline{u_i})- R(u_i, \overline{u_1}, u_1 - u_i, \overline{u_1} -  \overline{u_i}) \\
= {1 \over 2} \{ R(u_1 + u_i, \overline{u_1} + \overline{u_i}, u_1
+ u_i,\overline{u_1} + \overline{u_i} )
   -2 R(u_1 -u_i,\overline{u_1} - \overline{u_i}, u_1 + u_i,\overline{u_1} + \overline{u_i}  )
  \\
  \qquad \qquad \qquad \qquad  + R ((u_1 -u_i,\overline{u_1} - \overline{u_i}, u_1 - u_i,\overline{u_1} - \overline{u_i} ) \}
\\
 = {1 \over 2 }\cdot 4 \{R(A_i,\overline{A_i}, A_i,\overline{A_i}) -2 R(B_i, \overline{B_i}, A_i, \overline{A_i}) +
 R(B_i,\overline{B_i},B_i,\overline{B_i})\},
\end{array}
\]
where
\[
A_i = {1\over \sqrt{2}} (u_1 + u_i),\qquad {\rm and}\;\; B_i =
{1\over \sqrt{2}} (u_1 -u_i).
\]

If $i \neq 1$, then both $A_i$ and $B_i$ are unitary vectors in
$T^{1,0}M\;$ and $A_i \perp B_i.\;$ Comparing the last formula
with the formula for sectional curvature $K(e_1, J e_i),$ we
obtain
\[
K(e_1, J e_i) =  {1\over 4} \{R(A,\overline{A}, A,\overline{A}) -2
R(B, \overline{B}, A, \overline{A}) +
R(B,\overline{B},B,\overline{B})\}, \qquad \forall \;
i=1,2,\cdots, n.
\]
In particular when $i=2,$  we have (note that $A=A_2$ and $B=B_2$)
\[
K(e_1,J e_2) = K(w_1,w_2) =  {1\over 4} \{R(A,\overline{A},
A,\overline{A}) -2 R(B, \overline{B}, A, \overline{A}) +
R(B,\overline{B},B,\overline{B})\}.
\]
\end{proof}



\noindent Department of mathematics, Princeton University,
Princeton, NJ 08540, USA; \\
\noindent xiu@math.princeton.edu\\

\noindent Department of mathematics, MIT, Cambridge, MA
02139-4307, USA;\\
\noindent tian@math.mit.edu
\end{document}